\documentclass[review,3p,times,12pt]{elsarticle}
\usepackage{lineno,hyperref,amsmath,amsthm}
\usepackage[noend]{algpseudocode}
\usepackage{algorithmicx,algorithm}
\usepackage{booktabs}
\usepackage{color}
\usepackage{graphics}
\usepackage{subfig}
\usepackage{amsfonts,amssymb}
\usepackage{mathrsfs}
\usepackage{bm}

\renewcommand{\eqref}[1]{\textup{{\normalfont Eq.\,(\ref{#1}}\normalfont)}}
\allowdisplaybreaks[4]
\newtheorem{remark}{Remark}
\newcommand\figref{Fig.\,\ref}

\begin{document}

\begin{frontmatter}
    
    \title{A stochastic LATIN method for stochastic and parameterized elastoplastic analysis}

    \author[a,c]{Zhibao Zheng}

    \author[b,c]{David N{\'e}ron}
    
    \author[a,c]{Udo Nackenhorst}
    
    \address[a]{Leibniz Universit{\"{a}}t Hannover, Institute of Mechanics and Computational Mechanics, Appelstra{\ss}e 9a, 30167 Hannover, Germany}
    
    \address[b]{Universit{\'e} Paris-Saclay, CentraleSup{\'e}lec, ENS Paris-Saclay, CNRS, LMPS - Laboratoire de M{\'e}canique Paris-Saclay, 91190 Gif-sur-Yvette, France}

    \address[c]{International Research Training Group 2657 (IRTG 2657): Computational Mechanics Techniques in High Dimensions}

    \begin{abstract}
        The LATIN method has been developed and successfully applied to a variety of deterministic problems, but few work has been developed for nonlinear stochastic problems.
        This paper presents a \emph{stochastic LATIN method} to solve stochastic and/or parameterized elastoplastic problems.
        To this end, the stochastic solution is decoupled into spatial, temporal and stochastic spaces, and approximated by the sum of a set of products of triplets of spatial functions, temporal functions and random variables.
        Each triplet is then calculated in a greedy way using a stochastic LATIN iteration.
        The high efficiency of the proposed method relies on two aspects: The nonlinearity is efficiently handled by inheriting advantages of the classical LATIN method, and the randomness and/or parameters are effectively treated by a sample-based approximation of stochastic spaces.
        Further, the proposed method is not sensitive to the stochastic and/or parametric dimensions of inputs due to the sample description of stochastic spaces.
        It can thus be applied to high-dimensional stochastic and parameterized problems.
        Four numerical examples demonstrate the promising performance of the proposed stochastic LATIN method.
        
    \end{abstract}
    
    
\end{frontmatter}

\section{Introduction} \label{section:Intro}

	Developments in basic theories, numerical techniques and computing hardware have made it possible to solve high-resolution models in various science and engineering problems.
    However, in the presence of uncertainty in input data, such as random material properties, random external forces and random geometries \cite{stefanou2009stochastic, dai2023new, zheng2023simulation}, it may not always be desirable to augment the model complexity if the intrinsic predictive capability of the model is not satisfactory.
    The considerable influence of uncertainties on system behavior has led to the significant development of dedicated numerical techniques for uncertainty propagation \cite{le2010spectral, ghanem2017handbook}.
    Further, in many physical problems, model parameters are typically modeled as interval values instead of random variables/fields, known as parameterized problems \cite{hesthaven2016certified}.
    Generally speaking, parameterized problems are a simplified class of stochastic problems, where random inputs are given by uniform random variables (i.e. interval parameter values).
    In this paper, we focus on efficient and accurate numerical techniques for solving elastoplastic problems with random and parametric inputs.
    Since the parameter spaces of parameterized problems can be embedded into probability spaces, both stochastic and parameterized problems will be discussed in the context of a unified probability framework.

    Although deterministic elastoplastic problems have been well understood and solved effectively \cite{simo2006computational, de2012nonlinear}, the development of dedicated algorithms for stochastic elastoplasticity is still insufficient.
    Over the last two decades, several methods have been developed for this purpose.
    The Monte Carlo simulation (MCS) and its variants \cite{stefanou2009stochastic, papadrakakis1996robust, graham2011quasi} are general-purpose methods for various stochastic problems.
    This method only requires existing deterministic solvers and is therefore easy to be implemented.
    It can be applied to high-dimensional stochastic problems due to the independence from stochastic dimensions.
    However, a large number of deterministic simulations are required to achieve high-accuracy stochastic solutions, which is prohibitively expensive for large-scale stochastic elastoplastic problems.
    The well-known polynomial chaos (PC)-based spectral method \cite{ghanem2003stochastic, xiu2002wiener} has been extended to solve stochastic elastoplastic problems. 
    In PC-based methods, the stochastic solution is expanded over the PC basis and the stochastic Galerkin approach is then adopted to transform the original stochastic problem into an augmented deterministic equation whose size is much larger than the original stochastic problem.
    A direct extension of the PC-based method to stochastic elastoplastic problems can be found in \cite{rosic2011stochastic}.
    However, the classical PC expansion cannot well capture nonsmooth stochastic solutions of stochastic elastoplastic problems. 
    Several modifications have been developed to avoid this issue.  
    In \cite{anders1999stochastic, anders2001three}, stochastic elastoplastic constitutive equations are approximated by introducing two fictitious bounding and then solved by PC expansions.
    The Fokker-Planck-Kolmogorov equation based approach is proposed in \cite{sett2011stochastic, karapiperis2016fokker} to propagate uncertainty in stochastic elastoplastic constitutive models, and the PC expansion is used to approximate and solve unknown stochastic solutions.
    Stochastic elastoplastic problems are treated as stochastic variational inequalities in \cite{arnst2012variational}.
    The nonlinear convex programming coupled with PC expansions and stochastic collocation methods is used to solve stochastic solutions.
    The multi-element PC expansion is adopted to capture the nonsmooth behavior of stochastic elastoplastic problems in \cite{basmaji2022sepp}.
    Further, a stochastic Newton method is presented in \cite{zheng2023nonlinear} to solve (high-dimensional) stochastic elastoplastic problems.
    A time discretization framework is used to discretize the stochastic elastoplastic problem as (pseudo) time-dependent stochastic nonlinear problems.
    For each time step, a stochastic Newton iteration is proposed to linearize the stochastic nonlinear problem as a series of linear stochastic finite element equations, which are further solved by a weakly intrusive stochastic finite element method.
    There are also other attempts to solve stochastic elastoplastic problems, such as the stochastic perturbation method \cite{strkakowski2019stochastic} and the stochastic collocation methods \cite{dannert2022investigations, ajith2022novel}.

    We highlight that the computational burden of stochastic elastoplastic problems comes from the coupling of randomness and nonlinearity.
    However, all methods mentioned above focus on the effective treatment of randomness, while the nonlinearity is only handled by the Newton linearization.
    A natural idea is to reduce the computational effort by introducing more efficient numerical techniques to take into account for nonlinearities.
    The LArge Time INcrement (LATIN) method \cite{ladeveze2012nonlinear} is considered for this purpose in this paper.
    As a powerful nonlinear solver, the LATIN method has been applied to a variety of problems, e.g., large deformation problems \cite{bellenger2001phenomenological}, plastic and viscoplastic problems \cite{relun2013model, bergheau2016proper}, contact problems \cite{ribeaucourt2007new, boucard2011parallel}, damage problems \cite{bhattacharyya2018latin}, parameter identification \cite{allix2002new, nguyen2008robust}, multiscale and multiphysics problems \cite{neron2008computational, saavedra2017enhanced}, etc.
    Unlike the Newton method, which relies on time discretization and transforms history-dependent nonlinear problems into a sequence of nonlinear problems at each time step, the LATIN method decomposes the original nonlinear problem into two subproblems, called global and local stages. Alternating solutions of the two subproblems with specified search directions are used to solve a series of corrections of unknown solutions over both spatial and temporal spaces.
    Different choices of the subproblems and search directions induce different variants of LATIN methods \cite{scanff2021study}, but the two subproblems typically correspond to a linear global problem used to capture spatial solutions and a nonlinear local problem used to satisfy the constitutive evolution, which reduces the computational effort greatly compared to the Newton method due to the decoupling of nonlinear constitutive evolution and spatial solution.
    Extensions of the classical LATIN method to parameterized nonlinear problems have been extensively studied in \cite{vitse2014virtual, neron2015time,  courard2016integration}.
    These methods achieve accurate and efficient solutions via the neat treatment of parameterized inputs, but their extensions to stochastic problems need further investigation.

    In this paper, we extend the classical LATIN method to solve stochastic nonlinear problems, named the stochastic LATIN method, which solves both stochastic and parametric elastoplastic problems by a unified probabilistic framework.
    Similar to the classical case, the stochastic LATIN method decomposes the original stochastic elastoplastic problem into global and local stages, which are used to capture spatial solutions and perform stochastic elastoplastic constitutive evolution, respectively.    
    The search directions of both two stages can be inherited from the classical LATIN method.
    The unknown stochastic solution in both global and local stages is approximated by a randomized version of the proper generalization decomposition (PGD) \cite{ladeveze2012nonlinear, nouy2007generalized, chinesta2011short, chinesta2013proper}, that is, the stochastic solution is represented as the sum of a set of products of triplets of spatial functions, temporal functions and random variables.
    A greedy way is then used to sequentially calculate the triplets.
    For each triplet, the global and local stages correspond to solving global and linear deterministic problems for spatial components of the stochastic solution and local and nonlinear stochastic problems for temporal and random components of the stochastic solution, respectively.
    Different from existing similar PGD methods such as the PC expansion-based generalized spectral decomposition methods \cite{nouy2007generalized, nouy2009generalized}, we approximate the random components of stochastic solutions only using a small number of random sample realizations, which requires less computational effort and easily performs stochastic elastoplastic constitutive evolution.
    However, it is noted that the computational accuracy of above stochastic solution is less accurate since it is solved in a sequential way and only uses a few random samples.
    To avoid this issue, an update stage is further introduced to recalculate the random and temporal components of stochastic solutions using a reduced-order stochastic elastoplastic equation that is constructed by the obtained spatial functions.
    Furthermore, it will be shown in the studied examples that the sample approximation of stochastic spaces is not sensitive to the dimensions of random/parametric inputs.
    Thus, the proposed method can be applied to high-dimensional random and parameterized elastoplastic problems without modification.

    The paper is organized as follows: Section~\ref{section:StoEP} gives the basic setting and the key issues of stochastic and parameterized elastoplastic problems.
    The approximation of stochastic solutions is also presented in this section.
    Section~\ref{section:Sto_LATIN}
    introduces a stochastic LATIN method to solve stochastic elastoplastic problems.
    Computational implementations and algorithm flowchart of the proposed stochastic LATIN method are elaborated in Section~\ref{section:Scheme_Implementation}.
    Following that, four numerical examples are given to demonstrate the performance of the proposed method in Section~\ref{section:Examples},
    and conclusions and discussions are presented in Section~\ref{section:Conclusion}.

\section{Stochastic elastoplastic problems} \label{section:StoEP}

\subsection{Stochastic elastoplastic equation}

    In this paper, we solve elastoplastic equations with uncertain inputs and find the stochastic displacement ${\bm u}\left( {\bf x}, t, \theta \right): \Omega \times {\cal T} \times \Theta \rightarrow \mathbb{R}^d$ such that the following equation holds for $\forall {\bf x} \in \Omega$, $\forall t \in {\cal T}$ and $\forall \theta \in \Theta$,

    \vspace{-1.0em}
    \begin{equation} \label{eq:EP_probelm}
        \left \{  
        \begin{aligned}
            - \triangledown \cdot \left[ {\bm \sigma} \left( {\bm u}\left( {\bf x}, t, \theta \right), t, \theta \right) \right] &= {\bm f}\left( {{\bf{x}}, t, \theta } \right) ~~&{\rm{in}}~~~&\Omega \\
            {\bm \sigma} \left( {\bm u}\left( {\bf x}, t, \theta \right), t, \theta \right) \cdot {\bf{n}} &= {\bm h}\left( {\bf x}, t, \theta \right)~~&{\rm{on}}~~~&{\Gamma _N}\\
            {\bm u}\left( {\bf x}, t, \theta \right) &= {\overline{{\bm u}}}\left( {\bf x}, t, \theta \right)~~&{\rm{on}}~~~&{\Gamma _D}  
        \end{aligned}  
        \right. ,
    \end{equation}
    where the spatial dimension of the geometric domain $\Omega \subset \mathbb{R}^{d}$ may be $d=1,2,3$,
    and the (pseudo) time variable $t$ is related to history-dependent stochastic elastoplastic constitutive models.
    The random event $\theta \in \Theta$ is defined in a probability space $( \Theta, \Xi, \cal{P} )$, where $\Theta$ denotes the space of elementary events, $\Xi$ is a $\sigma$-algebra defined on $\Theta$ and $\cal{P}$ is the probability measure.
    Note that the parameters in parameterized elastoplastic problems are given by intervals rather than defined in a probability space.
    To handle both random and parametric inputs using a unified numerical framework, we consider the parameters as a set of uniform distribution random variables in this paper.
    Hereinafter, we slightly abuse the notation, and "random" and "$\theta$" are adopted to describe both random and parametric cases.
    The stochastic stress tensor ${\bm \sigma} \left( {\bm u}\left( {\bf x}, t, \theta \right), t, \theta \right) \in \mathbb{R}^{d \times d}$ is a nonlinear function of the stochastic strain tensor ${\bm \varepsilon} \left( {\bm u}\left( {\bf x}, t, \theta \right) \right) \in \mathbb{R}^{d \times d}$ and the stochastic displacement ${\bm u} \left( {\bf x}, t, \theta \right) \in \mathbb{R}^d$, which is also associated with the stochastic elastoplastic constitutive models.
    ${\bm f}\left( {\bf x}, t, \theta \right) \in \mathbb{R}^d$ is the (pseudo) time-dependent stochastic external force.
    ${\Gamma _N}$ and ${\Gamma _D}$ are boundary segments associated with the Neumann boundary condition ${\bm h}\left( {\bf x}, t, \theta \right)$ and the Dirichlet boundary condition $\overline{{\bm u}}\left( {\bf x}, t, \theta \right)$.
    Without loss of generality, only the homogeneous Dirichlet boundary condition ${\overline{\bm u}}\left( {\bf x}, t, \theta \right) = 0$ is considered in this paper.

    The stochastic elastoplastic problem under consideration is assumed to be geometrically linear
    and takes the following linear and small stochastic strain tensor
    
    \vspace{-1.0em}
    \begin{equation} \label{eq:sto_strain}
        {\bm \varepsilon} \left( {\bm u}\left( {\bf x}, t, \theta \right) \right) = \frac{1}{2}\left( {\nabla {\bm u}\left( {\bf x}, t, \theta \right) + \left( {\nabla}{\bm u}\left( {\bf x}, t, \theta \right) \right)^{\rm T} } \right),
    \end{equation}
    which can be further decomposed into two parts

    \vspace{-1.0em}
    \begin{equation} \label{eq:strain_ep}
        {\bm \varepsilon}\left( {\bm u}\left( {\bf x}, t, \theta \right) \right) = {{\bm \varepsilon} _{\rm e}}\left( {\bm u}\left( {\bf x}, t, \theta \right), \theta \right) + {{\bm \varepsilon} _{\rm p}}\left( {\bm u}\left( {\bf x}, t, \theta \right), \theta \right),
    \end{equation}
    where ${{\bm \varepsilon} _{\rm e}}\left( {\bm u}\left( {\bf x}, t, \theta \right) \right)$ and ${{\bm \varepsilon} _{\rm p}}\left( {\bm u}\left( {\bf x}, t, \theta \right) \right)$ denote the stochastic elastic strain tensor and the stochastic plastic strain tensor, respectively.
    Furthermore, the stochastic stress tensor ${\bm \sigma}\left( {\bm u}\left( {\bf x}, t, \theta \right), t, \theta \right)$ can be represented using the stochastic elastic strain tensor ${{\bm \varepsilon} _{\rm e}}\left( {\bm u}\left( {\bf x}, t, \theta \right) \right)$
    
    \vspace{-1.0em}
    \begin{equation} \label{eq:stress_strain}
        {\bm \sigma}\left( {{\bm u}\left( {\bf x}, t, \theta \right), t, \theta } \right) = {\bm C}\left( {\bf x}, \theta \right) {\bm \varepsilon} _{\rm e}\left( {\bm u}\left( {\bf x}, t, \theta \right), \theta \right) = {\bm C}\left( {\bf x}, \theta \right) \left( {\bm \varepsilon}\left( {\bm u}\left( {\bf x}, t, \theta \right) \right) - {\bm \varepsilon} _{\rm p}\left( {\bm u}\left( {\bf x}, t, \theta \right), \theta \right) \right),
    \end{equation}
    where ${\bm C}\left( {\bf x}, \theta \right)$ is the fourth-order elastic tensor and may depend on the spatial position ${\bf x}$.
    It has different forms in different background problems.
    Without loss of generality, the following form is adopted in this paper
    
    \vspace{-1.0em}
    \begin{equation} \label{eq:Ce}
        {\bm C}\left( {\bf x}, \theta \right) = \kappa\left( {\bf x}, \theta \right) {\bf I} \otimes {\bf I} + 2\mu\left( {\bf x}, \theta \right){\bf I}_{\rm dev},
    \end{equation}
    where the stochastic bulk modulus $\kappa\left( {\bf x}, \theta \right)$ and the stochastic shear modulus $\mu\left( {\bf x}, \theta \right)$ may be modeled as deterministic quantities, parameters, random variables, random fields, or a mixture of them. 
    The fourth-order tensor ${{\bf I}_{\rm dev}}$ is given by ${{\bf I}_{\rm dev}} = {{\bf I}_4} - \frac{1}{3}{\bf I} \otimes {\bf I}$, where the second-order identity tensor ${\bf I} \in \mathbb{R}^{d \times d}$, the fourth-order tensor ${{\bf I}_4} = {\delta _{ik}}{\delta _{jl}}{{\bf e}_i} \otimes {{\bf e}_j} \otimes {{\bf e}_k} \otimes {{\bf e}_l}$, and ${\delta _{ik}}$ is the Kronecker delta.

    Further, we adopt a stochastic elastoplastic constitutive relation obtained by introducing uncertainties into the classical von Mises yield criterion with the linear kinematic hardening \cite{de2012nonlinear, zheng2023nonlinear}.
    Nevertheless, more general elastoplastic constitutive relations with uncertainties can also be dealt with by the proposed method, such as the Drucker-Prager yield criterion and the isotropic hardening under uncertainties \cite{de2012nonlinear}.
    In the randomized von Mises yield and linear kinematic hardening criteria under consideration, the connection of the stochastic thermodynamical force ${\bm \beta}\left( {\bf x}, t, \theta \right)$ and the internal variable ${\bm \rho}\left( {\bf x}, t, \theta \right)$ is given by
    
    \vspace{-1.0em}
    \begin{equation} \label{eq:inter_pat}
        {\bm \beta}\left( {\bf x}, t, \theta \right) = \varpi\left( {\bf x}, t, \theta \right){\bm \rho}\left( {\bf x}, t, \theta \right),
    \end{equation}
    where the positive coefficient $\varpi\left( {\bf x}, t, \theta \right)$ can be modeled as a deterministic value, parameter, random variable, or random field in practice.
    Typically, we can formulate the evolution of stochastic strain and stochastic internal variables as

    \vspace{-1.0em}
    \begin{equation} \label{eq:evo}
        \left( {{\bm{\dot \varepsilon}} _{\rm p}}\left( {\bm u}\left( {\bf x}, t, \theta \right), t, \theta \right), {\bm{\dot \rho}}\left( {\bf x}, t, \theta \right) \right) = {\mathscr T} \left( {\bm \sigma}\left( {{\bm u}\left( {\bf x}, t, \theta \right), t,\theta } \right), {\bm \beta}\left( {\bf x}, t, \theta \right) \right),
    \end{equation}
    where ${\mathscr T}\left( \cdot \right)$ is the evolution operator.
    We only consider an abstract from here and a detailed stochastic elastoplastic constitutive evolution will be discussed in Section \ref{subsection:SEP_time}.
    Further, the stochastic yield function is given by 
    
    \vspace{-1.0em}
    \begin{equation} \label{eq:VM_yield}
        \Psi \left( {\bm \sigma}\left( {\bf x}, t, \theta \right) ,{\bm \beta}\left( {\bf x}, t, \theta \right), {\bf x}, \theta \right) = \left| {{\bf I}_{\rm dev}} {{\bm \sigma} \left( {\bf x}, t, \theta \right) - {\bm \beta} \left( {\bf x}, t, \theta \right)} \right| - {\sigma _{\rm Y}}\left( {\bf x}, \theta \right) ,
    \end{equation}
    where ${\sigma _{\rm Y}}\left( {\bf x}, \theta \right)$ is the deterministic/stochastic yield stress that may also be given by a deterministic value, parameter, random variable, or random field, and ${{\bf I}_{\rm dev}} {\bm \sigma} \left( {\bf x}, t, \theta \right)$ is considered as the stochastic stress deviator of the stochastic stress tensor ${\bm \sigma} \left( {\bf x}, t, \theta \right)$.

    In this way, it is seen that due to uncertainties in the external force, the material properties, the boundary conditions or the elastoplastic constitutive relation, the displacement, strain and stress are all stochastic quantities.
    An efficient and accurate stochastic solution of \eqref{eq:EP_probelm} remains challenging.
    On one hand, different from deterministic elastoplastic problems, a spatial point may be in elastic and plastic states simultaneously when we execute the judgment of \eqref{eq:VM_yield}.
    In other words, the spatial point is in elastic states for some sample realizations of $\theta$, but in plastic states for other sample realizations.
    It is not easy to capture the possible elastic or plastic states accurately.
    This issue becomes even more difficult when the random input has high dimensions and large variability.
    On another hand, the coupling of nonlinearity and randomness makes \eqref{eq:EP_probelm} computationally demanding.
    Efficient solution algorithms cannot be easily extended from existing deterministic nonlinear solvers.
    The development of more advanced numerical solvers is thus expected.

\subsection{Weak form of stochastic elastoplastic equations}
    To solve \eqref{eq:EP_probelm}, let us consider its weak form in the spatial domain written as
    
    \vspace{-1.0em}
    \begin{equation} \label{eq:weak_form}
        {\mathscr W} \left( {\bm u}\left( {\bf x}, t, \theta \right),{\bm v}\left( {\bf{x}} \right), \theta \right) = {\mathscr F} \left( {\bm v}\left( {\bf{x}} \right), t, \theta \right), \quad \forall {\bm v}\left( {\bf{x}} \right) \in {\cal V},
    \end{equation}
    where ${\bm v}\left( {\bf x} \right)$ is the test function and ${\cal V}$ is the Hilbert space with smooth displacement vectors vanishing on the boundary ${\Gamma _D}$.
    The stochastic terms ${\mathscr W} \left( {\bm u}\left( {\bf x}, t, \theta \right),{\bm v}\left( {\bf{x}} \right), \theta \right)$ and ${\mathscr F} \left( {\bm v}\left( {\bf{x}} \right), t, \theta \right)$ are respectively given by

    \vspace{-1.0em}
    \begin{align}
        {\mathscr W} \left( {\bm u}\left( {\bf x}, t, \theta \right),{\bm v}\left( {\bf{x}} \right), \theta \right) &= \int_\Omega  {{\bm \sigma} \left( {{\bm u}\left( {\bf x}, t, \theta \right),t,\theta } \right) \colon {\bm \varepsilon} \left( {{\bm v}\left( {\bf{x}} \right)} \right){\rm d}{\bf{x}}}, \label{eq:weak_W} \\
        {\mathscr F} \left( {\bm v}\left( {\bf{x}} \right), t, \theta \right) &= \int_\Omega {\bm f}\left( {\bf x}, t, \theta \right) \cdot {\bm v}\left( {\bf x} \right){\rm d}{\bf x}  + \int_{\Gamma _N} {\bm h}\left( {\bf x}, t, \theta \right) \cdot {\bm v}\left( {\bf x} \right){\rm d}{\bf s}. \label{eq:weak_F}
    \end{align}
    \eqref{eq:weak_form} is a (pseudo) time-dependent nonlinear stochastic problem.
    Note that it is just a nominal weak form and does not involve detailed stochastic elastoplastic constitutive evolution. 
    Different techniques for linearization and iteration yield different approaches.
    For instance, in our previous work \cite{zheng2023nonlinear}, a stochastic Newton iteration based nonlinear stochastic finite element method is proposed by a time discretization framework coupled with linearized stochastic stresses. 
    However, its computational effort strongly depends on the time discretization framework.
    A large number of linearized stochastic finite element equations need to be solved if many time steps are involved.
    In order to reduce the computational burden of the stochastic Newton iteration, a stochastic LATIN method is developed in this paper.

\subsection{Stochastic solution approximation}
    To solve \eqref{eq:weak_form}, we first approximate the stochastic solution ${\bm u}\left( {\bf x}, t, \theta \right)$ using the following $k$-term PGD with random variable coefficients

    \vspace{-1.0em}
    \begin{equation} \label{eq:u_app}
        {\bm u}_k\left( {\bf x}, t,\theta \right) = \sum\limits_{l=1}^k \lambda_l\left( \theta \right) g_l\left( t \right) {\bm d}_l\left( {\bf x} \right) = {\bm u}_{k-1}\left( {\bf x}, t,\theta \right) + \lambda_k\left( \theta \right) g_k\left( t \right) {\bm d}_k\left( {\bf x} \right),
    \end{equation}
    where ${\bm u}_{k-1}\left( {\bf x}, t,\theta \right)$ is the $\left( k-1 \right)$-term approximation of the stochastic solution ${\bm u}\left( {\bf x}, t, \theta \right)$, $\left\{ \lambda_l\left( \theta \right) \right\}_{l=1}^k$ are scalar random variables, $\left\{ g_l\left( t \right) \right\}_{l=1}^k$ are temporal functions and $\left\{ {\bm d}_l\left( {\bf x} \right) \right\}_{l=1}^k$ are spatial functions.
    \eqref{eq:u_app} decouples the stochastic solution into spatial, temporal and stochastic spaces.
    All components are possibly determined in their individual spaces.
    However, all triplets $\left\{ \lambda_l\left( \theta \right), g_l\left( t \right), {\bm d}_l\left( {\bf x} \right) \right\}_{l=1}^k$ are not known a priori in practice.
    Following the PGD spirit, a greedy procedure is adopted to calculate the triplet $\left\{ \lambda_l\left( \theta \right), g_l\left( t \right), {\bm d}_l\left( {\bf x} \right) \right\}$, $l = 1, 2, \cdots$ one by one in a sequential way.
    Specifically, we assume that the $\left( k -1 \right)$-term approximation ${\bm u}_{k-1}\left( {\bf x}, t,\theta \right)$ in \eqref{eq:u_app} has been known and only the $k$-th triplet $\left\{ \lambda_k\left( \theta \right), g_k\left( t \right), {\bm d}_k\left( {\bf x} \right) \right\}$ is unknown.
    They are then solved in an iterative way.
    Following this idea, an investigation of solving stochastic elastoplastic problems was carried out within the framework of stochastic Newton iteration \cite{zheng2023nonlinear}.
    In the present paper, the triplets will be determined iteratively in a stochastic LATIN framework.

    Further, parameterized and randomized PGD methods have been developed to solve parameterized and stochastic nonlinear problems \cite{neron2015time, nouy2009generalized, vitse2019dealing}.
    The increment of \eqref{eq:u_app} is concentrated on the effective treatment of randomness.
    On one hand, the randomness can be captured using only a set of scalar random variables, which provides a pathway to propagate and quantify complex uncertain inputs in a simple way.
    On another hand, we adopt random sample realizations to describe the random variables $\left\{ \lambda_l\left( \theta \right) \right\}_{l=1}^k$ in the numerical implementations.
    The effectiveness of this strategy has been verified in \cite{zheng2023nonlinear, zheng2023stochastic}.
    Compared to explicit approximations of these random variables such as the PC expansion \cite{nouy2009generalized}, the proposed method can examine the elastic or plastic state of a spatial point via \eqref{eq:VM_yield} in a way similar to deterministic elastoplastic problems.
    Since the random sample realizations are not sensitive to the stochastic dimensions and can be used to approximate high-dimensional random inputs without modification.
    Thus, the curse of dimensionality arising in high-dimensional stochastic and parametric problems can be avoided successfully.
    Further, the proposed method is weakly intrusive, that is, \eqref{eq:u_app} is an intrusive approximation, but the implementation of the random variables $\left\{ \lambda_l\left( \theta \right) \right\}_{l=1}^k$ is non-intrusive.
    Existing numerical techniques and codes can be reused to a great extent, e.g., (large-scale) linear equation solvers, finite element discretization and assembly, etc.

\section{A stochastic LATIN framework for nonlinear stochastic analysis} \label{section:Sto_LATIN}
\begin{figure}[ht]
    \centering
    \includegraphics[width=0.7\linewidth]{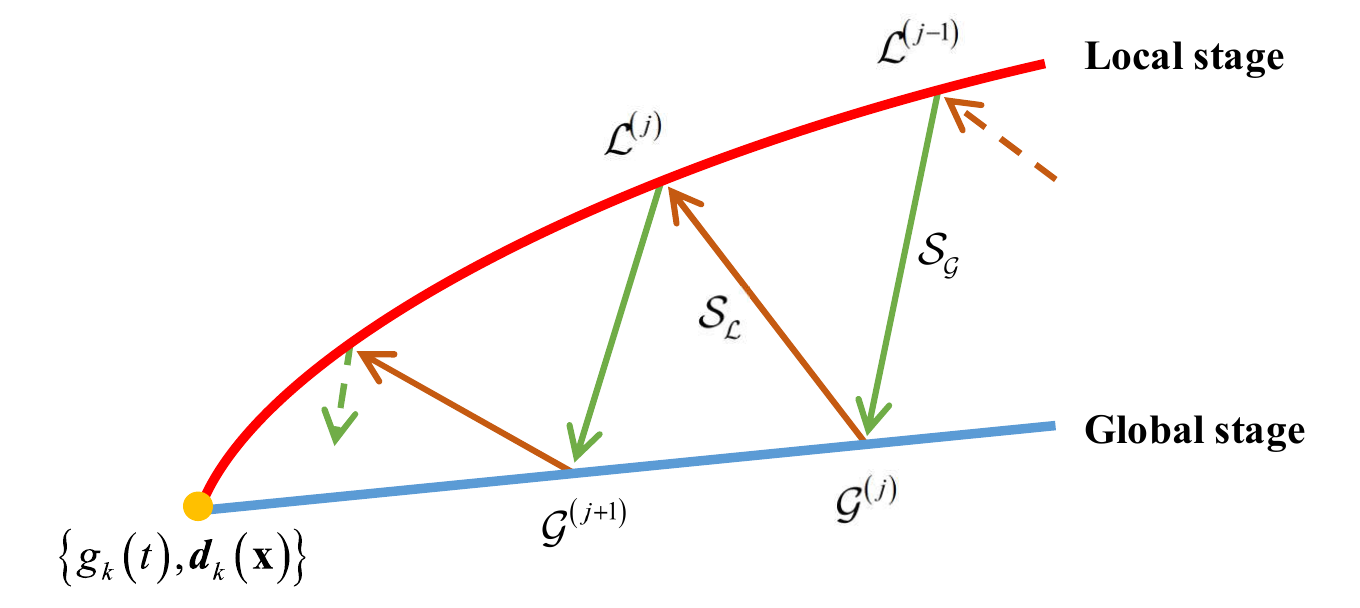}
    \caption{Deterministic LATIN iteration for the $k$-th triplet $\left\{ g_k\left( t \right), {\bm d}_k\left( {\bf x} \right) \right\}$.}
    \label{fig_Der_LATIN_iteration}
\end{figure}

    In this section, we focus on extending the classical LATIN method \cite{ladeveze2012nonlinear} to a stochastic LATIN iteration, which is used herein to efficiently treat the coupling of nonlinearity and randomness.
    Let us simply recall the idea of classical LATIN iteration, which is usually combined with the deterministic PGD expansion.
    For the $k$-th pair $\left\{ g_k\left( t \right), {\bm d}_k\left( {\bf x} \right) \right\}$ used to approximate the solution of deterministic problems, the LATIN iteration decomposes the original problem into a local and a global stage, as shown in \figref{fig_Der_LATIN_iteration}.
    Two search directions ${\cal S}_{\cal L}$ and ${\cal S}_{\cal G}$ are used to solve the solutions of local and global stages, respectively.
    Possible choices of the search directions can be found in \cite{scanff2021study}.
    Typically, the global stage is only related to global linear spatial equations,
    and the local stage is associated with local nonlinear equations and used to satisfy the elastoplastic constitutive evolution.
    Also, it is noted that the PGD approximation is only used at the global stage to reduce the computational effort.
    It is not used in the local stage since the calculation of the local stage is cheap enough.
    We can obtain the solution of the pair $\left\{ g_k\left( t \right), {\bm d}_k\left( {\bf x} \right) \right\}$ after several alternating iterations of local and global stages.
    Compared to the Newton linearization, the LATIN method decouples spatial solutions and nonlinear constitutive evolution. 
    It allows solving nonlinear problems with less global linearized equations in spatial domains and local nonlinear equations in temporal domains, which thus involves less computational effort.
    More theoretical details of the classical LATIN method such as convergence analysis can be found in \cite{ladeveze2012nonlinear}.

\begin{figure}[ht]
    \centering
    \includegraphics[width=0.7\linewidth]{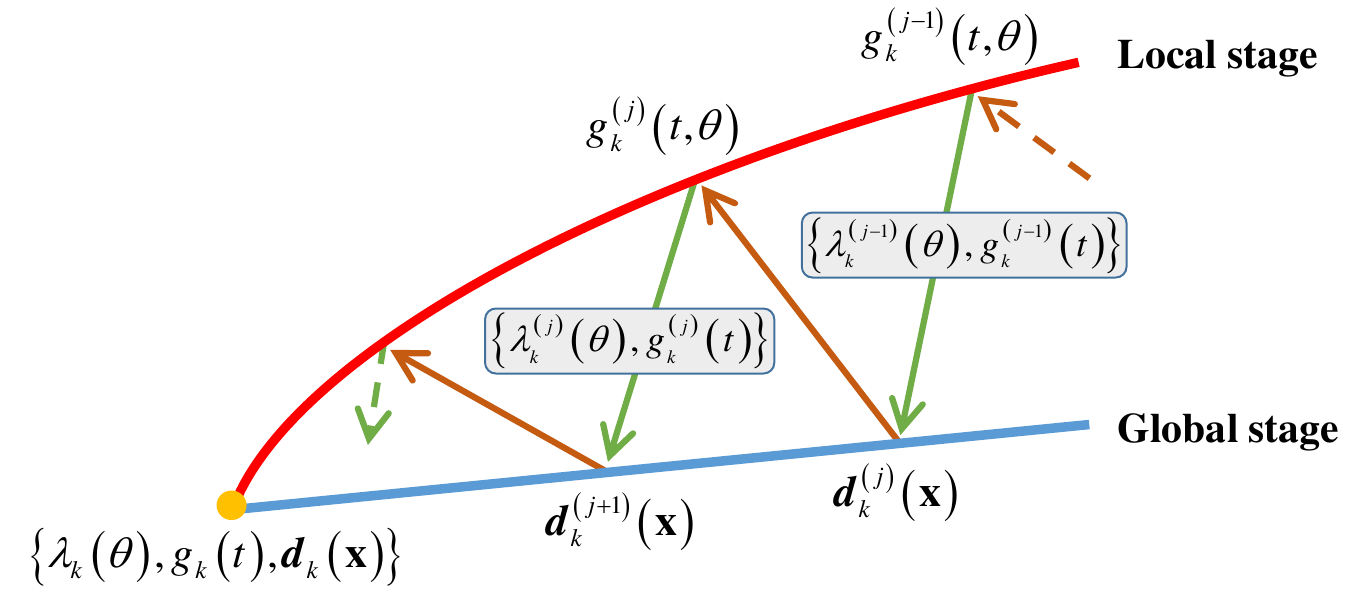}
    \caption{Stochastic LATIN iteration framework for the $k$-th triplet $\left\{ \lambda_k\left( \theta \right), g_k\left( t \right), {\bm d}_k\left( {\bf x} \right) \right\}$.}
    \label{fig_LATIN_iteration}
\end{figure}

    Analogous to the above deterministic LATIN method, a stochastic LATIN iterative process is proposed to determine the $k$-th triplet component $\left\{ \lambda_k\left( \theta \right), g_k\left( t \right), {\bm d}_k\left( {\bf x} \right) \right\}$ in \eqref{eq:u_app}.
    As illustrated in \figref{fig_LATIN_iteration}, the stochastic LATIN iteration decomposes the original stochastic problem into two stages, including a linear and deterministic global stage and a nonlinear and stochastic local stage, where the global stage is used to calculate the spatial function, and the local stage is used to satisfy the stochastic elastoplastic constitutive evolution and determine the random and temporal components.
    In both stages, additional concerns about the effective treatment of random inputs are involved.
    Moreover, since expensive computation effort is required at the local stage to handle the randomness, we also adopt the PGD at the local stage to reduce the effort.
    In this way, the PGD is used to both global and local stages, which is different from the classical LATIN method where the PGD is only used at the global stage.
    In this section, we only provide a basic stochastic LATIN framework.
    The implementation details will be shown in Section \ref{section:Scheme_Implementation}.

\subsection{Global stage of stochastic LATIN iteration} \label{subsection:Global_stage}
    In the global stage of the proposed stochastic LATIN iteration, we iteratively calculate the spatial function ${\bm d}_k\left( {\bf x} \right)$ of the $k$-th triplet $\left\{ \lambda_k\left( \theta \right), g_k\left( t \right), {\bm d}_k\left( {\bf x} \right) \right\}$ via linear and deterministic global problems.
    To this end, we assume that the pair $\left\{ \lambda_k^{\left( j-1 \right)}\left( \theta \right), g_k^{\left( j-1 \right)} \left( t \right) \right\}$ has been known and then solve ${\bm d}_k^{\left( j \right)}\left( {\bf x} \right)$.
    The search direction ${\cal S}_{\cal G}$ of the global stage involves random and/or parametric inputs, and is set to find a linear stochastic stress increment on the stochastic stress ${\bm \sigma}\left( {\bm \varepsilon} \left( {\bm u}_{k-1}\left( {\bf x}, t, \theta \right) \right), t, \theta \right)$ obtained by the $\left( k-1 \right)$-term approximation of the stochastic solution.
    Therefore, the satisfaction of stochastic elastoplastic constitutive evolution is not requested at this stage.
    The original weak form (\ref{eq:weak_form}) is further reformulated as a deterministic linear problem under the given random variable $\lambda_k^{\left( j-1 \right)}\left( \theta \right)$ and temporal function $g_k^{\left( j-1 \right)} \left( t \right)$,
    
    \vspace{-1.0em}
    \begin{equation} \label{eq:SL_d}
        {\mathscr W}_{k,{\bm d}}^{\left( j \right)} \left( {\bm d}_k^{\left( j \right)}\left( {\bf x} \right),{\bm v}\left( {\bf{x}} \right) \right) = {\mathscr F}_{k,{\bm d}}^{\left( j \right)} \left( {\bm v}\left( {\bf x} \right) \right) ,
    \end{equation}
    where the deterministic terms ${\mathscr W}_{k,{\bm d}}^{\left( j \right)} \left( {\bm d}_k^{\left( j \right)}\left( {\bf x} \right),{\bm v}\left( {\bf{x}} \right) \right) \in \mathbb{R}$ and ${\mathscr F}_{k,{\bm d}}^{\left( j \right)} \left( {\bm v}\left( {\bf x} \right) \right) \in \mathbb{R}$ are given by adopting the classical Galerkin process in the temporal space and the stochastic Galerkin process in the stochastic space \cite{ghanem2003stochastic, xiu2002wiener}

    \vspace{-1.0em}
    \begin{align}
        {\mathscr W}_{k,{\bm d}}^{\left( j \right)} \left( {\bm d}_k^{\left( j \right)}\left( {\bf x} \right),{\bm v}\left( {\bf{x}} \right) \right) &= \iiint_{{\Omega} \times {\cal T} \times \Theta} \left( \lambda_k^{\left( j-1 \right)}\left( \theta \right) g_k^{\left( j-1 \right)} \left( t \right) \right)^2 \left[ {\bm C}\left( {\bf x}, \theta \right) {\bm \varepsilon} \left( {\bm d}_k^{\left( j \right)}\left( {\bf x} \right) \right) \colon {\bm \varepsilon} \left( {\bm v}\left( {\bf x} \right) \right) \right] {\rm d}{\bf x} {\rm d}t {\rm d}{\cal P}\left( \theta \right), \label{eq:SL_d_W} \\
        {\mathscr F}_{k,{\bm d}}^{\left( j \right)} \left( {\bm v}\left( {\bf x} \right) \right)&= \iint_{{\cal T} \times \Theta} \lambda_k^{\left( j-1 \right)}\left( \theta \right) g_k^{\left( j-1 \right)} \left( t \right) {\mathscr F} \left( {\bm v}\left( {\bf x} \right), t, \theta \right) {\rm d}t {\rm d}{\cal P}\left( \theta \right) \nonumber \\
        &\hspace{-0.95cm} - \iiint_{{\Omega} \times {\cal T} \times \Theta} \lambda_k^{\left( j-1 \right)}\left( \theta \right) g_k^{\left( j-1 \right)} \left( t \right) {\bm \sigma}\left( {\bm \varepsilon} \left( {\bm u}_{k-1}\left( {\bf x}, t,\theta \right) \right), t, \theta \right) \colon {\bm \varepsilon} \left( {\bm v}\left( {\bf x} \right) \right){\rm d}{\bf x} {\rm d}t {\rm d}{\cal P}\left( \theta \right), \label{eq:SL_d_F} 
    \end{align}
    where the stochastic term $\lambda_k^{\left( j-1 \right)}\left( \theta \right) g_k^{\left( j-1 \right)} \left( t \right) {\bm C}\left( {\bf x}, \theta \right) {\bm \varepsilon} \left( {\bm d}_k^{\left( j \right)}\left( {\bf x} \right) \right)$ is the linear stochastic stress increment on the stochastic stress ${\bm \sigma}\left( {\bm \varepsilon} \left( {\bm u}_{k-1}\left( {\bf x}, t, \theta \right) \right), t, \theta \right)$.
    ${\bm d}_k^{\left( j \right)}\left( {\bf x} \right)$ represents the $j$-th intermediate iteration of ${\bm d}_k\left( {\bf x} \right)$, which is the unknown solution to be solved in this stage.
    The known $\lambda_k^{\left( j-1 \right)}\left( \theta \right)$ and $g_k^{\left( j-1 \right)} \left( t \right)$ are the $\left( j-1 \right)$-th intermediate iterations of $\lambda_k\left( \theta \right)$ and $g_k\left( t \right)$ obtained by the previous iterative step $j-1$ (or given initial values).
    In this way, \eqref{eq:SL_d} is a deterministic system of linear equations about ${\bm d}_k^{\left( j \right)}\left( {\bf x} \right)$, which can be discretized and solved by means of existing numerical techniques.
    It is much more efficient and has less storage memory compared with solving a nonlinear problem of ${\bm d}_k^{\left( j \right)}\left( {\bf x} \right)$.
    Detailed numerical implementations of the global stage will be presented in Section \ref{subsec:Im_glo}.

\subsection{Local stage of stochastic LATIN iteration}\label{subsection:Local_stage}
    On the basis of the spatial function ${\bm d}_k^{\left( j \right)}\left( {\bf x} \right)$ solved by \eqref{eq:SL_d} in the global stage, the pair $\left\{ \lambda_k^{\left( j \right)}\left( \theta \right), g_k^{\left( j \right)}\left( t \right) \right\}$ is updated using a nonlinear and stochastic local problem in the local stage, and the search direction ${\cal S}_{\cal L}$ is set to perform the stochastic elastoplastic constitutive evolution.
    As shown in \figref{fig_LATIN_iteration}, we first solve a coupled temporal-stochastic solution $\hbar_k^{\left( j \right)}\left( t,\theta \right)$ instead of the temporal-stochastic separation form $\left\{ \lambda_k^{\left( j \right)}\left( \theta \right), g_k^{\left( j \right)} \left( t \right) \right\}$. 
    To this end, the weak form (\ref{eq:weak_form}) is reformulated as a local nonlinear stochastic problem based on the known spatial function ${\bm d}_k^{\left( j \right)}\left( {\bf x} \right)$ 

    \vspace{-1.0em}
    \begin{equation} \label{eq:SL_g}
        {\mathscr W}_{k,\hbar}^{\left( j \right)} \left( \hbar_k^{\left( j \right)}\left( t,\theta \right), \theta \right) = {\mathscr F}_{k,\hbar}^{\left( j \right)} \left( t,\theta \right),
    \end{equation}
    where the stochastic terms ${\mathscr W}_{k,\hbar}^{\left( j \right)} \left( \hbar_k^{\left( j \right)}\left( t,\theta \right), \theta \right)$ and ${\mathscr F}_{k,\hbar}^{\left( j \right)} \left( t,\theta \right)$ are given by adopting the classical Galerkin process in the spatial domain
    
    \vspace{-1.0em}
    \begin{align}
        {\mathscr W}_{k,\hbar}^{\left( j \right)} \left( \hbar_k^{\left( j \right)}\left( t,\theta \right), \theta \right) &= \int_\Omega {\bm \sigma} \left( {\bm u}_{k-1}\left( {\bf x}, t, \theta \right) + \hbar_k^{\left( j \right)}\left( t,\theta \right) {\bm d}_k^{\left( j \right)}\left( {\bf x} \right), t, \theta \right) \colon {\bm \varepsilon} \left( {\bm d}_k^{\left( j \right)}\left( {\bf x} \right) \right){\rm d}{\bf x}, \label{eq:SL_g_W} \\
        {\mathscr F}_{k,\hbar}^{\left( j \right)} \left( t,\theta \right) &= \int_\Omega  {\bm f}\left( {\bf x}, t, \theta \right) \cdot {\bm d}_k^{\left( j \right)}\left( {\bf x} \right){\rm d}{\bf x}  + \int_{\Gamma_N} {\bm h}\left( {\bf x}, t, \theta \right) \cdot {\bm d}_k^{\left( j \right)}\left( {\bf x} \right){\rm d}{\bf s}, \label{eq:SL_g_F}
    \end{align}
    where the stochastic stress tensor ${\bm \sigma} \left( {\bm u}_{k-1}\left( {\bf x}, t, \theta \right) + \hbar_k^{\left( j \right)}\left( t,\theta \right) {\bm d}_k^{\left( j \right)}\left( {\bf x} \right), t, \theta \right)$ is calculated via the nonlinear stochastic elastoplastic constitutive evolution.
    In this way, \eqref{eq:SL_g} is a one-dimensional and time-dependent nonlinear stochastic problem about $\hbar_k^{\left( j \right)}\left( t,\theta \right)$.
    It can be solved cheaply by classical nonlinear numerical solvers, such as the fixed-point iteration and the Newton iteration.    
    Further, to achieve a separated pair $\lambda_k^{\left( j \right)}\left( \theta \right)$ and $g_k^{\left( j \right)}\left( t \right)$, we decompose $\hbar_k^{\left( j \right)}\left( t, \theta \right)$ as follows

    \vspace{-1.0em}
    \begin{equation} \label{eq:SL_g_lam_decouple}
        \hbar_k^{\left( j \right)}\left( t, \theta \right) = \lambda_k^{\left( j \right)}\left( \theta \right) g_k^{\left( j \right)}\left( t \right),
    \end{equation}
    which greatly saves storage memory, but only provides a rough approximation of $\hbar_k^{\left( j \right)}\left( t, \theta \right)$.
    In our experience, the decomposition \eqref{eq:SL_g_lam_decouple} has less influence on the accuracy of the final stochastic solution ${\bm u} \left( {\bf x}, t, \theta \right)$.
    In fact, the approximation accuracy of the stochastic solution ${\bm u}_k \left( {\bf x}, t, \theta \right)$ will be improved as the retained number $k$ increases although the pair $\left\{ \lambda_k^{\left( j \right)}\left( \theta \right), g_k^{\left( j \right)}\left( t \right) \right\}$ is less accurate for the current triplet compared to $\hbar_k^{\left( j \right)}\left( t, \theta \right)$.
    Further, we will show that \eqref{eq:SL_g_lam_decouple} can be considered as a rank-1 randomized singular value decomposition (SVD) in Section \ref{subsec:non_it_local_stage}.

    Repeatedly solving the problems (\ref{eq:SL_d}), (\ref{eq:SL_g}) and (\ref{eq:SL_g_lam_decouple}) until a specified precision is reached, we can determine the converged solution of the $k$-th triplet $\left\{ \lambda_k\left( \theta \right), g_k\left( t \right), {\bm d}_k\left( {\bf x} \right) \right\}$.
    The next triplets $\left\{ \lambda_{k+1}\left( \theta \right), g_{k+1}\left( t \right), {\bm d}_{k+1}\left( {\bf x} \right) \right\}$, $\cdots$ can be solved in a similar way.
    The stochastic solution ${\bm u}\left( {\bf x}, t, \theta \right)$ can be accurately approximated by retaining enough triplets.
    Further, when considering each random sample realization $\left\{ \lambda_l\left( \theta^* \right) \right\}_{l=1}^k$, $\theta^* \in \Theta$ individually, \eqref{eq:u_app} degenerates into the deterministic PGD approximation.
    Thus, convergence and accuracy of the stochastic solution approximation can be inherited from the classical LATIN iteration in a probability sense.

\begin{remark}
    We remark that in \eqref{eq:SL_g}, we can also solve the decoupled pair $\left\{ \lambda_k^{\left( j \right)}\left( \theta \right), g_k^{\left( j \right)}\left( t \right) \right\}$ directly, which simply requires reformulating the left side of \eqref{eq:SL_g} as
    
    \vspace{-1.0em}
    \begin{equation} \label{eq:SL_g_W_2}
        {\mathscr W}_{k,\hbar}^{\left( j \right)} \left( \lambda_k^{\left( j \right)}\left( \theta \right), g_k^{\left( j \right)}\left( t \right), \theta \right) = \int_\Omega {\bm \sigma} \left( {\bm u}_{k-1}\left( {\bf x}, t, \theta \right) + \lambda_k^{\left( j \right)}\left( \theta \right) g_k^{\left( j \right)}\left( t \right) {\bm d}_k^{\left( j \right)}\left( {\bf x} \right),t,\theta \right) \colon {\bm \varepsilon} \left( {\bm d}_k^{\left( j \right)}\left( {\bf x} \right) \right){\rm d}{\bf x}, 
    \end{equation}
    which is a nonlinear stochastic problem about both $\lambda_k^{\left( j \right)}\left( \theta \right)$ and $g_k^{\left( j \right)}\left( t \right)$.
    Compared to \eqref{eq:SL_g}, \eqref{eq:SL_g_W_2} is computationally expensive since extra effort is required to solve a time-dependent nonlinear problem about $g_k^{\left( j \right)}\left( t \right)$ and a stochastic nonlinear problem about $\lambda_k^{\left( j \right)}\left( \theta \right)$.
\end{remark}

\subsection{Update stage of stochastic solutions} \label{subsection:Update_stage}
    It is noted that the above iteration solves the triplets $\left\{ \lambda_l\left( \theta \right), g_l\left( t \right), {\bm d}_l\left( {\bf x} \right) \right\}_{l=1}^k$ in a sequential way.
    In our experience, this iteration cannot achieve a very accurate approximation of the stochastic solution for all possible random sample realizations.
    To improve the accuracy of the final stochastic solution, we consider ${\bm D}\left( {\bf x} \right) = \left[ {\bm d}_1\left( {\bf x} \right), \cdots, {\bm d}_k\left( {\bf x} \right) \right]$ as a set of reduced basis functions and approximate the stochastic solution ${\bm u}_k\left( {\bf x}, t,\theta \right)$ as
    
    \vspace{-1.0em}
    \begin{equation} \label{eq:u_rom}
        {\bm u}_k\left( {\bf x}, t,\theta \right) = {\bm D}\left( {\bf x} \right) {\bm g}\left( t, \theta \right),
    \end{equation}
    where the temporal-stochastic coupling functions ${\bm g}\left( t, \theta \right) = \left[ g_1\left( t, \theta \right), \cdots, g_k\left( t, \theta \right) \right]^{\rm T} \in \mathbb{R}^k$ are unknown.
    Then, we recalculate ${\bm g}\left( t, \theta \right)$ by solving the following reduced-order nonlinear stochastic elastoplastic problem

    \vspace{-1.0em}
    \begin{equation} \label{eq:weak_form_rom}
        {\mathscr W}_r \left( {\bm g}\left( t, \theta \right), {{\bm v}\left( {\bf{x}} \right)}, \theta \right) = {\mathscr F}_r \left( {\bm v}\left( {\bf x} \right), \theta \right),
    \end{equation}
    where the stochastic terms ${\mathscr W}_r \left( {\bm g}\left( t, \theta \right), {{\bm v}\left( {\bf{x}} \right)}, \theta \right)$ and ${\mathscr F}_r \left( {\bm v}\left( {\bf x} \right), \theta \right)$ are given by

    \vspace{-1.0em}
    \begin{align}
        {\mathscr W}_r \left( {\bm g}\left( t, \theta \right), {\bm v}\left( {\bf{x}} \right), \theta \right) &= \int_\Omega  {{\bm \sigma} \left( {\bm D}\left( {\bf x} \right) {\bm g}\left( t, \theta \right), t, \theta \right) \colon {\bm \varepsilon} \left( {{\bm v}\left( {\bf{x}} \right)} \right){\rm d}{\bf{x}}}, \label{eq:weak_W_rom} \\
        {\mathscr F}_r \left( {\bm v}\left( {\bf{x}} \right), \theta \right) &= \int_\Omega {\bm f}\left( {\bf x}, t, \theta \right) \cdot {\bm v}\left( {\bf x} \right){\rm d}{\bf x}  + \int_{\Gamma _N} {\bm h}\left( {\bf x}, t, \theta \right) \cdot {\bm v}\left( {\bf x} \right){\rm d}{\bf s}. \label{eq:weak_F_rom}
    \end{align}
    It is noted that ${\bm D}\left( {\bf x} \right)$ is also used as the basis functions of approximate space of the test function ${\bm v}\left( {\bf x} \right)$, that is, ${\bm v}\left( {\bf x} \right) = \sum_{i=1}^k {\bm d}_i\left( {\bf x} \right) v_i$ and $\left\{ v_i \right\}_{i=1}^k$ are the unknown coefficients.
    The solution of \eqref{eq:weak_form_rom} is cheap since its size $k$ is much smaller than the original stochastic problem (\ref{eq:weak_form}).

\section{Computational scheme implementation} \label{section:Scheme_Implementation}
    In this section, we present the numerical implementation details for each stage of the above stochastic LATIN iteration.
    The evolution of stochastic elastoplastic constitutive model is elaborated in the implementation of local stage.
    The algorithm flowchart is also provided to clearly present the proposed stochastic LATIN method.

\subsection{Implementation of global stage} \label{subsec:Im_glo}
    To solve \eqref{eq:SL_d}, we adopt the finite element method for spatial discretization.
    Other discretization techniques are also available for the purpose, such as the finite volume method, the isogeometric approximation and the virtual element method \cite{cottrell2009isogeometric, antonietti2022virtual}. 
    Finite element approximations of the spatial functions ${\bm d}_k^{\left( j \right)}\left( {\bf x} \right)$ and ${\bm v}\left( {\bf x} \right)$ are given by
    
    \vspace{-1.0em}
    \begin{equation} \label{eq:dv_fe}
        {\bm d}_k^{\left( j \right)}\left( {\bf x} \right) = \sum\limits_{i=1}^n {\bm \varphi}_i\left( {\bf x} \right) d_{k,i}^{\left( j \right)} = {\bf \Phi}\left( {\bf x} \right) {\bf d}_k^{\left( j \right)}, \quad
        {\bm v}\left( {\bf x} \right) = \sum\limits_{i=1}^n {\bm \varphi}_i\left( {\bf x} \right) v_i = {\bf \Phi}\left( {\bf x} \right) {\bf v},
    \end{equation}
    where $n$ is the number of spatial degrees of freedom, 
    ${\bf \Phi}\left( {\bf x} \right) = \left[ {\bm \varphi}_1\left( {\bf x} \right), \cdots, {\bm \varphi}_n\left( {\bf x} \right) \right]$ is a set of finite element basis functions, and ${\bf d}_k^{\left( j \right)} = \left[ d_{k,1}^{\left( j \right)}, \cdots, d_{k,n}^{\left( j \right)} \right]^{\rm T} \in \mathbb{R}^n$ and ${\bf v} = \left[ v_1, \cdots, v_n \right]^{\rm T} \in \mathbb{R}^n$ are unknown vectors of the discretized nodal values.
    Substituting \eqref{eq:dv_fe} into Eqs.~(\ref{eq:SL_d_W}) and (\ref{eq:SL_d_F}), \eqref{eq:SL_d} is transformed into a deterministic system of linear equations
    
    \vspace{-1.0em}
    \begin{equation} \label{eq:KdF}
        {\bf K}_k^{\left( j \right)} {\bf d}_k^{\left( j \right)} = {\bf F}_k^{\left( j \right)},
    \end{equation}
    where the elements of the deterministic matrix ${\bf K}_k^{\left( j \right)} \in \mathbb{R}^{n \times n}$ are given by

    \vspace{-1.0em}
    \begin{equation} \label{eq:KdF_K}
        {\bf K}_{k,sr}^{\left( j \right)} = \int_{\Omega} {\bm C}_{k}^{\left( j \right)}\left( {\bf x} \right) {\bm \varepsilon} \left( {\bm \varphi}_r\left( {\bf x} \right) \right) \colon {\bm \varepsilon} \left( {\bm \varphi}_s\left( {\bf x} \right) \right) {\rm d}{\bf x} , \quad s,r = 1, \cdots, n,
    \end{equation}
    and the deterministic tensor ${\bm C}_{k}^{\left( j \right)}$ is given by
    
    \vspace{-1.0em}
    \begin{equation} \label{eq:KdF_C}
        {\bm C}_{k}^{\left( j \right)}\left( {\bf x} \right) = \iint_{{\cal T} \times \Theta} \left( \lambda_k^{\left( j-1 \right)}\left( \theta \right) g_k^{\left( j-1 \right)} \left( t \right) \right)^2 {\bm C}\left( {\bf x}, \theta \right) {\rm d}t {\rm d}{\cal P}\left( \theta \right),
    \end{equation}
    which involves the integration on both temporal and stochastic domains.
    In the numerical implementation, the time variable $t$ is discretized as $n_t$ time steps, i.e. $t \in \left[t_1, \cdots, t_{n_t} \right]$.
    We adopt the rectangle formula for the temporal integration and a sample-based approach for the stochastic integration, which corresponds to

    \vspace{-1.0em}
    \begin{equation} \label{eq:KdF_C_Com}
        {\bm C}_{k}^{\left( j \right)}\left( {\bf x} \right) = {\widehat{\mathbb{E}}}\left\{ \left( \lambda_k^{\left( j-1 \right)} \left( {\widehat{\bm \theta}} \right) \right)^2 \odot {\bm C}\left( {\bf x}, {\widehat{\bm \theta}} \right) \right\} \sum\limits_{i=1}^{n_t} \left( g_k^{\left( j-1 \right)} \left( t_i^* \right) \right)^2 \Delta t_i,
    \end{equation}
    where $\Delta t_i$ is the time interval and $t_i^* = \frac{t_i + t_{i+1}}{2}$ is the middle point of $t_i$ and $t_{i+1}$.
    $\lambda_k^{\left( j-1 \right)} \left( {\widehat{\bm \theta}} \right)$ and ${\bm C}\left( {\bf x}, {\widehat{\bm \theta}} \right)$ are a set of sample realizations of the random variable $\lambda_k^{\left( j-1 \right)} \left( \theta \right)$ and the random tensor ${\bm C}\left( {\bf x}, \theta \right)$, respectively.
    $\odot$ represents the element-wise multiplication operator, and ${\widehat{\mathbb{E}}}\left\{ \cdot \right\}$ is the expectation operator of sample realizations, with ${\widehat{\mathbb{E}}}\left\{ \left( \lambda_k^{\left( j-1 \right)} \left( {\widehat{\bm \theta}} \right) \right)^2 \odot {\bm C}\left( {\bf x}, {\widehat{\bm \theta}} \right) \right\} = \frac{1}{n_s}\sum_{i=1}^{n_s} \left( \lambda_k^{\left( j-1 \right)} \left( \theta^{\left( i \right)} \right) \right)^2 {\bm C}\left( {\bf x}, \theta^{\left( i \right)} \right)$, where $n_s$ is the number of random sample realizations.  
    Similarly, the elements of the deterministic force vector ${\bf F}_k^{\left( j \right)} \in \mathbb{R}^n$ in \eqref{eq:KdF} are calculated by

    \vspace{-1.0em}
    \begin{align} \label{eq:KdF_F}
        {\bf F}_{k,s}^{\left( j \right)} = &\int_{\Omega} \left[  \sum\limits_{i=1}^{n_t} g_k^{\left( j-1 \right)} \left( t_i^* \right) {\widehat{\mathbb{E}}}\left\{ \lambda_k^{\left( j-1 \right)} \left( {\widehat{\bm \theta}} \right) \odot {\bm f}\left( {\bf x}, t_i^*, {\widehat{\bm \theta}} \right) \right\} \Delta t_i \right] \cdot {\bm \varphi}_s\left( {\bf x} \right){\rm d}{\bf x} \nonumber \\
        &+ \int_{\Gamma _N} \left[  \sum\limits_{i=1}^{n_t} g_k^{\left( j-1 \right)} \left( t_i^* \right) {\widehat{\mathbb{E}}}\left\{ \lambda_k^{\left( j-1 \right)} \left( {\widehat{\bm \theta}} \right) \odot {\bm h}\left( {\bf x}, t_i^*, {\widehat{\bm \theta}} \right) \right\} \Delta t_i \right] \cdot {\bm \varphi}_s\left( {\bf x} \right){\rm d}{\bf s} \nonumber \\
        & - \int_\Omega \left[ \sum\limits_{i=1}^{n_t} g_k^{\left( j-1 \right)} \left( t_i^* \right){\widehat{\mathbb{E}}}\left\{ \lambda_k^{\left( j-1 \right)} \left( {\widehat{\bm \theta}} \right) \odot {\bm \sigma}\left( {\bm \varepsilon} \left( {\bm u}_{k-1}\left( {\bf x}, t_i^*, {\widehat{\bm \theta}} \right) \right), t_i^*, {\widehat{\bm \theta}} \right) \right\} \Delta t_i \right] \colon {\bm \varepsilon} \left( {\bm \varphi}_s\left( {\bf x} \right) \right){\rm d}{\bf x} 
    \end{align}
    for $s = 1, \cdots, n$.

    It is noted that the above numerical integration is less accurate in both temporal and stochastic spaces, but has little influence on the accuracy of the final stochastic solution ${\bm u}\left( {\bf x}, t, \theta \right)$.
    For the stochastic integration, we do not need to accurately calculate the random variable $\lambda_k^{\left( j-1 \right)} \left( \theta \right)$,
    and only the expectation ${\widehat{\mathbb{E}}}\left\{ \left( \lambda_k^{\left( j-1 \right)} \left( {\widehat{\bm \theta}} \right) \right)^2 \odot {\bm C}\left( {\bf x}, {\widehat{\bm \theta}} \right) \right\}$ (and other expectations in \eqref{eq:KdF_F}) is expected to be accurate.
    By using the iterative process proposed in the previous section, the expectation value can be gradually approximated to a good accuracy although the random variable $\lambda_k^{\left( j-1 \right)} \left( \theta \right)$ may not be determined accurately.
    Also, for both temporal and stochastic integration, the approximation accuracy increases as the retained number $k$ of triplets increases.
    Further, we can finally improve the accuracy of the stochastic solution by solving the reduced-order stochastic elastoplastic problem (\ref{eq:weak_form_rom}).

    Furthermore, the random tensor ${\bm C}\left( {\bf x}, \theta \right)$ has an affine form in many cases.
    If the affine form does not exist, we can achieve an \eqref{eq:C_expansion}-like approximation by using the methods for simulating random fields, such as the Karhunen-Lo{\`e}ve expansion and the PC expansion \cite{zheng2017simulation, sakamoto2002polynomial, zheng2021sample}.
    In this paper, we only consider the following affine form to explain the proposed method

    \vspace{-1.0em}
    \begin{equation} \label{eq:C_expansion}
        {\bm C}\left( {\bf x}, \theta \right) = \sum\limits_{i=1}^m \xi_i\left( \theta \right) {\bm C}_i\left( {\bf x} \right) ,
    \end{equation}
    where $\left\{ \xi_i\left( \theta \right) \right\}_{i=1}^m$ are random variables and $\left\{ {\bm C}_i\left( {\bf x} \right) \right\}_{i=1}^m$ are deterministic tensor components.
    In this way, \eqref{eq:KdF_C_Com} can be reformulated as

    \vspace{-1.0em}
    \begin{equation} \label{eq:KdF_C_ex}
        {\bm C}_{k}^{\left( j \right)}\left( {\bf x} \right) = \left( \sum\limits_{i=1}^m  \underbrace{{\widehat{\mathbb{E}}}\left\{ \left( \lambda_k^{\left( j-1 \right)} \left( {\widehat{\bm \theta}} \right) \right)^2 \odot \xi_i\left( {\widehat{\bm \theta}} \right) \right\}}_{=z_{k,i}^{\left( j \right)}} {\bm C}_i\left( {\bf x} \right) \right) \left( \underbrace{\sum\limits_{i=1}^{n_t} \left( g_k^{\left( j-1 \right)} \left( t_i^* \right) \right)^2 \Delta t_i}_{=1}  \right),
    \end{equation}
    where the deterministic scalar coefficients $z_{k,i}^{\left( j \right)} = {\widehat{\mathbb{E}}}\left\{ \left( \lambda_k^{\left( j-1 \right)} \left( {\widehat{\bm \theta}} \right) \right)^2 \odot \xi_i\left( {\widehat{\bm \theta}} \right) \right\}$ is introduced, and the normalization $\sum_{i=1}^{n_t} \left( g_k^{\left( j-1 \right)} \left( t_i^* \right) \right)^2 \Delta t_i = 1$ is prescribed.
    Thus, \eqref{eq:KdF} is rewritten as

    \vspace{-1.0em}
    \begin{equation} \label{eq:KdF_ex}
        \left( \sum\limits_{i=1}^m z_{k,i}^{\left( j \right)} {\bf K}_i \right) {\bf d}_k^{\left( j \right)} = {\bf F}_k^{\left( j \right)},
    \end{equation}
    where the deterministic matrices $\left\{ {\bf K}_i \in \mathbb{R}^{n \times n} \right\}_{i=1}^m$ are given by

    \vspace{-1.0em}
    \begin{equation} \label{eq:K_exp}
        {\bf K}_{i,sr} = \int_{\Omega} {\bm C}_i\left( {\bf x} \right) {\bm \varepsilon} \left( {\bm \varphi}_r\left( {\bf x} \right) \right) \colon {\bm \varepsilon} \left( {\bm \varphi}_s\left( {\bf x} \right) \right) {\rm d}{\bf x} , \quad s,r = 1, \cdots, n,
    \end{equation}
    which is only calculated once for all possible iterations, thus saving a lot of computational effort.
    During the iterative process, we only need to update the scalar coefficients $\left\{ z_{k,i}^{\left( j \right)} \right\}_{i=1}^m$ via the expectation operation of random sample vectors, whose cost is almost negligible even for very high-dimensional random inputs.

\subsection{Implementation of local stage} 
    In this section, we present the numerical implementation details for the local stage, including iterative treatment of the local nonlinear stochastic equation and the evolution of stochastic elastoplastic constitutive model.
    
\subsubsection{Nonlinear iteration of local stage} \label{subsec:non_it_local_stage}
    We further solve \eqref{eq:SL_g} based on the solution vector ${\bf d}_k^{\left( j \right)}$ solved by \eqref{eq:KdF}.
    The Newton linearization approach is adopted to deal with the nonlinearity and perform the evolution of time steps.
    To this end, we consider the following one-dimensional linearized stochastic equation of \eqref{eq:SL_g} at the time step $t_i$
    
    \vspace{-1.0em}
    \begin{equation} \label{eq:g_nl}
        a_k^{\left( j, m \right)}\left( t_i, \theta \right) \Delta g_k^{\left( j,m \right)}\left( t_i, \theta \right) = b_k^{\left( j, m \right)}\left( t_i, \theta \right),
    \end{equation}
    where $\Delta g_k^{\left( j,m \right)}\left( t_i, \theta \right) \in \mathbb{R}$ is the stochastic increment on the scalar stochastic solution $g_k^{\left( j,m \right)}\left( t_i, \theta \right)$ of the $m$-th Newton iteration.
    $a_k^{\left( j, m \right)}\left( t_i, \theta \right)$ and $b_k^{\left( j, m \right)}\left( t_i, \theta \right)$ are scalar random variables and given by the linearization of ${\mathscr W}_{k,\hbar}^{\left( j \right)} \left( \hbar_k^{\left( j \right)}\left( t,\theta \right), \theta \right)$ in \eqref{eq:SL_g}

    \vspace{-1.0em}
    \begin{align}
        a_k^{\left( j, m \right)}\left( t_i, \theta \right) &= \int_\Omega {\bm C}_{\rm T} \left( {\bm u}_{k-1}\left( {\bf x}, t_i, \theta \right) + g_k^{\left( j,m \right)}\left( t_i,\theta \right) {\bm d}_k^{\left( j \right)}\left( {\bf x} \right),\theta \right) {\bm \varepsilon} \left( {\bm d}_k^{\left( j \right)}\left( {\bf x} \right) \right) \colon {\bm \varepsilon} \left( {\bm d}_k^{\left( j \right)}\left( {\bf x} \right) \right){\rm d}{\bf x}, \label{eq:g_nl_a}\\
        b_k^{\left( j, m \right)}\left( t_i, \theta \right) &= {\mathscr F}_{k,\hbar}^{\left( j \right)} \left( t_i,\theta \right) - \int_\Omega {\bm \sigma} \left( {\bm u}_{k-1}\left( {\bf x}, t_i, \theta \right) + g_k^{\left( j,m \right)}\left( t_i,\theta \right) {\bm d}_k^{\left( j \right)}\left( {\bf x} \right), t_i, \theta \right) \colon {\bm \varepsilon} \left( {\bm d}_k^{\left( j \right)}\left( {\bf x} \right) \right){\rm d}{\bf x}, \label{eq:g_nl_b}
    \end{align}
    where ${\bm C}_{\rm T} = {{\partial {\bm{\sigma }}} \mathord{\left/ {\vphantom {{\partial {\bm{\sigma }}} {\partial {\bm{\varepsilon }}}}} \right. \kern-\nulldelimiterspace} {\partial {\bm{\varepsilon }}}}$ represents the tangent tensor used to linearize the nonlinear stochastic stress, which is related to the elastoplastic constitutive model and the time steps.
    An explicit representation of ${\bm C}_{\rm T}$ can be found in \eqref{eq:C_T}.
    Classical approaches for solving \eqref{eq:g_nl}, such as the PC method and the stochastic collocation method, are computationally complex and require dedicated algorithms to accurately and efficiently solve stochastic solutions.
    Most of these methods are also sensitive to the dimensions of random inputs and thus suffer from the curse of dimensionality for high-dimensional stochastic elastoplastic problems.
    To avoid these issues, the following sample-based non-intrusive method is used to calculate the random sample realizations of $\Delta g_k^{\left( j,m \right)}\left( t_i, \theta \right)$

    \vspace{-1.0em}
    \begin{equation} \label{eq:g_sol_sam}
        \Delta g_k^{\left( j,m \right)}\left( t_i, {\widehat{\bm \theta}} \right) = b_k^{\left( m,j \right)}\left( t_i, {\widehat{\bm \theta}} \right) \oslash a_k^{\left( m,j \right)}\left( t_i, {\widehat{\bm \theta}} \right) \in \mathbb{R}^{n_s},
    \end{equation}
    where $\oslash$ denotes the element-wise division of two vectors.
    $a_k^{\left( m \right)}\left( t_i, {\widehat{\bm \theta}} \right) \in \mathbb{R}^{n_s}$ and $b_k^{\left( m \right)}\left( t_i, {\widehat{\bm \theta}} \right) \in \mathbb{R}^{n_s}$ represent the sample realization vectors of the random variables $a_k^{\left( j, m \right)}\left( t_i, \theta \right)$ and $b_k^{\left( j, m \right)}\left( t_i, \theta \right)$, which is calculated using Eqs.~(\ref{eq:g_nl_a}) and (\ref{eq:g_nl_b}) for all possible sample realizations.
    Further, we update the stochastic solution $g_k^{\left( j,m+1 \right)}\left( t_i, \theta \right)$ with the sample realizations

    \vspace{-1.0em}
    \begin{equation} \label{eq:g_update}
        g_k^{\left( j,m+1 \right)}\left( t_i, {\widehat{\bm \theta}} \right) = g_k^{\left( j,m \right)}\left( t_i, {\widehat{\bm \theta}} \right) + \Delta g_k^{\left( j,m \right)}\left( t_i, {\widehat{\bm \theta}} \right) \in \mathbb{R}^{n_s}.
    \end{equation} 
    
    Repeatedly solving \eqref{eq:g_sol_sam} and \eqref{eq:g_update} until a specified precision is met, we can obtain the converged stochastic solution $g_k^{\left( j \right)}\left( t_i, {\widehat{\bm \theta}} \right) \in \mathbb{R}^{n_s}$.
    Looping all time steps we can obtain the sample realization matrix $g_k^{\left( j \right)}\left( t, {\widehat{\bm \theta}} \right) = \left[g_k^{\left( j \right)}\left( t_1, {\widehat{\bm \theta}} \right), \cdots, g_k^{\left( j,m+1 \right)}\left( t_{n_t}, {\widehat{\bm \theta}} \right) \right] \in \mathbb{R}^{n_s \times n_t}$ of the stochastic solution $\hbar_k^{\left( j \right)}\left( t, \theta \right)$.
    On this basis, \eqref{eq:SL_g_lam_decouple} is achieved via the rank-1 deterministic SVD of the sample realization matrix $g_k^{\left( j \right)}\left( t, {\widehat{\bm \theta}} \right)$

    \vspace{-1.0em}
    \begin{equation} \label{eq:g_lam_sam}
        g_k^{\left( j \right)}\left( t, {\widehat{\bm \theta}} \right) = \lambda_k^{\left( j \right)}\left( {\widehat{\bm \theta}} \right) g_k^{\left( j \right)}\left( t \right),
    \end{equation} 
    where $\lambda_k^{\left( j \right)}\left( {\widehat{\bm \theta}} \right) \in \mathbb{R}^{n_s}$ is the sample vector of the random variable $\lambda_k^{\left( j \right)}\left( \theta \right)$
    and $g_k^{\left( j \right)}\left( t \right) \in \mathbb{R}^{1 \times n_t}$ is the discretized vector of the temporal function.
    In a probability sense, \eqref{eq:SL_g_lam_decouple} is a rank-1 randomized SVD.
    It can be further extended to rank-$q$ randomized SVD, with $\hbar_k^{\left( j \right)}\left( t, \theta \right) = \sum_{i=1}^q \lambda_{k,i}^{\left( j \right)}\left( \theta \right) g_{k,i}^{\left( j \right)}\left( t \right)$, which can also be implemented with an \eqref{eq:g_lam_sam}-like rank-$q$ deterministic SVD.
    
\subsubsection{Elastoplastic constitutive evolution} \label{subsection:SEP_time}
    In this section, we perform temporal evolution of the stochastic elastoplastic constitutive model, which is associated with the linearized random coefficients $a_k^{\left( j, m \right)}\left( t_i, \theta \right)$ and $b_k^{\left( j, m \right)}\left( t_i, \theta \right)$ in \eqref{eq:g_nl}. 
    For each time step $t_i$, the stochastic solution ${\bm u}\left( {\bf x}, t_i, \theta \right)$ is updated by solving the local nonlinear stochastic problem (\ref{eq:g_nl}) that is dependent on the quantities $\left\{{\bm u}\left( {\bf x}, t_{i-1}, \theta \right), {\bm \varepsilon} _{{\rm p}}\left( {\bf x}, t_{i-1}, \theta \right), {\bm \beta} \left( {\bf x}, t_{i-1}, \theta \right) \right\}$ obtained from the previous time step $t_{i-1}$.
    For the stochastic elastoplastic constitutive model under consideration, an explicit evolution of the stochastic stress is given by 
    
    \vspace{-1.0em}
    \begin{equation} \label{eq:G_stress_strin}
        {{\bm \sigma}}\left({{\bm \varepsilon}}\left( {\bf x}, t_i, \theta \right); {\bm \varepsilon} _{{\rm p}}\left( {\bf x}, t_{i-1}, \theta \right), {\bm \beta}\left( {\bf x}, t_{i-1}, \theta \right) \right) = {{\bm \sigma} _{0}}\left( {\bf x}, t_i, \theta \right) - \frac{{2\mu \left( {\bf x}, \theta \right)}}{{2\mu \left( {\bf x}, \theta \right) + \varpi\left( {\bf x}, \theta \right)}}\left\langle {{\gamma}\left( {\bf x}, t_i, \theta \right)} \right\rangle {\bm \sigma} ^*\left( {\bf x}, t_i, \theta \right),
    \end{equation}
    where the operator $\left\langle \Box \right\rangle  = \max \left( {\Box, 0} \right)$ represents the maximum value of $\Box$ and 0, and $\gamma \left( {\bf x}, t_i, \theta \right) = 1 - \frac{{{\sigma _{\rm Y}}\left( {\bf x}, \theta \right)}}{{\left| {{\bm \sigma} ^*\left( {\bf x}, t_i, \theta \right)} \right|}}$.
    The stochastic stresses ${{\bm \sigma} _{0}}\left( {\bf x}, t_i, \theta \right)$ and ${\bm \sigma}^*\left( {\bf x}, t_i, \theta \right)$ are given by
    
    \vspace{-1.0em}
    \begin{align} 
        {{\bm \sigma} _{0}}\left( {\bf x}, t_i, \theta \right) &= {\bm C}\left( {\bf x}, \theta \right)\left( {{{\bm \varepsilon}}\left( {\bf x}, t_i, \theta \right) - {{\bm \varepsilon} _{{\rm p}}}\left( {\bf x}, t_{i-1}, \theta \right)} \right), \label{eq:Sto_stress_i} \\
        {\bm \sigma}^*\left( {\bf x}, t_i, \theta \right) &= {{\bf I}_{\rm dev}}{{\bm \sigma} _{0}}\left( {\bf x}, t_i, \theta \right) - {{\bm \beta}}\left( {\bf x}, t_{i-1}, \theta \right).\label{eq:Sto_stress_i_normal}
    \end{align}
    The tangent tensor ${\bm C}_{\rm T}$ in \eqref{eq:g_nl_a} is further given by \cite{zheng2023nonlinear, carstensen2002elastoviscoplastic}
    
    \vspace{-1.0em}
    \begin{align}
        &{\bm C}_{\rm T} \left({{\bm \varepsilon}}\left( {\bf x}, t_i, \theta \right); {\bm \varepsilon} _{{\rm p}}\left( {\bf x}, t_{i-1}, \theta \right), {\bm \beta}\left( {\bf x}, t_{i-1}, \theta \right) \right) = \frac{\partial {{\bm \sigma}}\left({{\bm \varepsilon}}\left( {\bf x}, t_i, \theta \right); {\bm \varepsilon} _{{\rm p}}\left( {\bf x}, t_{i-1}, \theta \right), {\bm \beta}\left( {\bf x}, t_{i-1}, \theta \right) \right)}{\partial {{\bm \varepsilon}}\left( {\bf x}, t_i, \theta \right)} \label{eq:C_T_0} \\
        = \; &{\bm C}\left( {\bf x}, \theta \right) - \frac{{4\mu {{\left( {\bf x}, \theta \right)}^2}}}{{2\mu \left( {\bf x}, \theta \right) + \varpi\left( {\bf x}, \theta \right)}}\left[ {{{\bf I}_{\rm dev}} +  \left( {\frac{1}{{{\gamma}\left( {\bf x}, t_i, \theta \right)} + \delta_{ {\gamma}\left( {\bf x}, t_i, \theta \right)} } - 1} \right) \frac{{{\bm \sigma}^*\left( {\bf x}, t_i, \theta \right) \otimes {\bm \sigma}^*\left( {\bf x}, t_i, \theta \right)}}{{{{\left| {{\bm \sigma}^*\left( {\bf x}, t_i, \theta \right)} \right|}^2}}}} \right]\left\langle {{\gamma}\left( {\bf x}, t_i, \theta \right)} \right\rangle, \label{eq:C_T}
    \end{align}
    where $\delta_{ {\gamma}\left( {\bf x}, t_i, \theta \right)}$ is the Kronecker delta meeting $\delta_{ {\gamma}\left( {\bf x}, t_i, \theta \right)} = 1$ for ${\gamma _i}\left( \theta \right) = 0$ and 0 for others.
    When ${\left| {{\bm \sigma}^*\left( {\bf x}, t_i, \theta \right)} \right|} \le {\sigma _{\rm Y}}\left( {\bf x}, \theta \right)$ (i.e. the true stress is less than the yield stress), ${\bm C}_{\rm T} = {\bm C}\left( {\bf x}, \theta \right)$ holds and we do not need to update ${\bm C}_{\rm T}$.
    Otherwise, the update of ${\bm C}_{\rm T}$ is necessary to perform the stochastic elastoplastic constitutive evolution.
    Since \eqref{eq:g_sol_sam} is solved using random samples, the proposed method can accurately execute the judgment ${\left| {{\bm \sigma}^*\left( {\bf x}, t_i, \theta^{\star} \right)} \right|} > {\sigma _{\rm Y}}\left( {\bf x}, \theta^{\star} \right)$ (and ${\left| {{\bm \sigma}^*\left( {\bf x}, t_i, \theta^{\star} \right)} \right|} \le {\sigma _{\rm Y}}\left( {\bf x}, \theta^{\star} \right)$) for each sample realization $\theta^{\star} \in \Theta$.
    Therefore, a way similar to deterministic elastoplastic constitutive evolution can still be adopted in the proposed method although random inputs are involved in stochastic elastoplastic problems.
    
    Further, ${{\bm \beta}}\left( {\bf x}, t_i, \theta \right)$ and ${{\bm \varepsilon} _{{\rm p}}}\left( {\bf x}, t_i, \theta \right)$ are then updated by using ${\bm \sigma}^*\left( {\bf x}, t_i, \theta \right)$ from current time step $t_i$ and ${{\bm \beta}}\left( {\bf x}, t_{i-1}, \theta \right)$, ${{\bm \varepsilon} _{{\rm p}}}\left( {\bf x}, t_{i-1}, \theta \right)$ from previous time step $t_{i-1}$
    
    \vspace{-1.0em}
    \begin{align}
        {{\bm \beta}}\left( {\bf x}, t_i, \theta \right) &= {{\bm \beta}}\left( {\bf x}, t_{i-1}, \theta \right) + \frac{\varpi\left( {\bf x}, \theta \right)}{{2\mu \left( {\bf x}, \theta \right) + \varpi\left( {\bf x}, \theta \right)}}\left\langle {{\gamma}\left( {\bf x}, t_i, \theta \right)} \right\rangle {\bm \sigma}^*\left( {\bf x}, t_i, \theta \right), \label{eq:beta_i} \\
        {{\bm \varepsilon} _{{\rm p}}}\left( {\bf x}, t_i, \theta \right) &= {{\bm \varepsilon} _{{\rm p}}}\left( {\bf x}, t_{i-1}, \theta \right) + \frac{1}{{2\mu \left( {\bf x}, \theta \right) + \varpi\left( {\bf x}, \theta \right)}}\left\langle {{\gamma}\left( {\bf x}, t_i, \theta \right)} \right\rangle {\bm \sigma}^*\left( {\bf x}, t_i, \theta \right). \label{eq:strin_i}
    \end{align}
    Similarly, the second terms of Eqs.~(\ref{eq:beta_i}) and (\ref{eq:strin_i}) are nonzero when ${\left| {{\bm \sigma}^*\left( {\bf x}, t_i, \theta \right)} \right|} > {\sigma _{\rm Y}}\left( {\bf x}, \theta \right)$ holds, the corresponding nodes are in plastic states. Otherwise, these second terms are equal to zero and we do not need to update Eqs.~(\ref{eq:beta_i}) and (\ref{eq:strin_i}).

\subsection{Implementation of update stage} 
    In this section, the Newton iteration is used to solve the reduced-order stochastic elastoplastic problem (\ref{eq:weak_form_rom}), which is cheap enough due to the small number $k$ of the reduced basis functions.
    In the practical implementation, we solve the problem (\ref{eq:weak_form_rom}) for each sample realization $\theta^{\left( q \right)}$ $q = 1, \cdots, n_s$.
    In this way, the total computational cost for $n_s$ random sample realizations is proportional to the number $n_s$.
    For each sample realization $\theta^{\left( q \right)}$ and the time step $t_i$, the following linearized equation of \eqref{eq:weak_form_rom} is generated by taking advantage of the Newton linearization and the classical Galerkin procedure

    \vspace{-1.0em}
    \begin{equation} \label{eq:kgf_rom}
        {\bf k}_{\rm T}\left( {\bf g}^{\left( m \right)}\left( t_i, \theta^{\left( q \right)} \right), \theta^{\left( q \right)} \right) \Delta {\bf g}^{\left( m \right)}\left( t_i, \theta^{\left( q \right)} \right) = {\bf f}_{\rm T}\left( {\bf g}^{\left( m \right)}\left( t_i, \theta^{\left( q \right)} \right), \theta^{\left( q \right)} \right),
    \end{equation} 
    where the reduced-order tangent matrix ${\bf k}_{{\rm T},sr}\left( {\bf g}^{\left( m \right)}\left( t_i, \theta^{\left( q \right)} \right), \theta^{\left( q \right)} \right) \in \mathbb{R}^{k \times k}$ 
    and the reduced-order vector ${\bf f}_{\rm T}\left( {\bf g}^{\left( m \right)}\left( t_i, \theta^{\left( q \right)} \right), \theta^{\left( q \right)} \right) \in \mathbb{R}^k$ are given by

    \vspace{-1.0em}
    \begin{align}
        {\bf k}_{{\rm T},sr}\left( {\bf g}^{\left( m \right)}\left( t_i, \theta^{\left( q \right)} \right), \theta^{\left( q \right)} \right) &= \int_\Omega {\bm C}_{\rm T} \left( {\bm D}\left( {\bf x} \right) {\bf g}^{\left( m \right)}\left( t_i, \theta^{\left( q \right)} \right), \theta^{\left( q \right)} \right) {\bm \varepsilon} \left( {\bm d}_r\left( {\bf x} \right) \right) \colon {\bm \varepsilon} \left( {\bm d}_s\left( {\bf x} \right) \right){\rm d}{\bf x}, \label{eq:k_T_rom} \\
        {\bf f}_{{\rm T},s}\left( {\bf g}^{\left( m \right)}\left( t_i, \theta^{\left( q \right)} \right), \theta^{\left( q \right)} \right) &= \int_\Omega  {\bm f}\left( {\bf x}, t_i, \theta^{\left( q \right)} \right) \cdot {\bm d}_s\left( {\bf x} \right){\rm d}{\bf x} + \int_{\Gamma_N} h\left( {\bf x}, t_i, \theta^{\left( q \right)} \right) \cdot {\bm d}_s\left( {\bf x} \right){\rm d}{\bf s} \nonumber \\
        &\qquad\qquad\qquad - \int_\Omega {\bm \sigma} \left( {\bm D}\left( {\bf x} \right) {\bf g}^{\left( m \right)}\left( t_i, \theta^{\left( q \right)} \right), t_i, \theta^{\left( q \right)} \right) \colon {\bm \varepsilon} \left( {\bm d}_s\left( {\bf x} \right) \right){\rm d}{\bf x}. \label{eq:f_rom}
    \end{align}
    The solution ${\bf g}\left( t_i, \theta^{\left( q \right)} \right)$ is then updated as follows until convergence
 
    \vspace{-1.0em}
    \begin{equation} \label{eq:g_rom_update}
        {\bf g}^{\left( m+1 \right)}\left( t_i, \theta^{\left( q \right)} \right) = {\bf g}^{\left( m \right)}\left( t_i, \theta^{\left( q \right)} \right) + \Delta {\bf g}^{\left( m \right)}\left( t_i, \theta^{\left( q \right)} \right) \in \mathbb{R}^k.
    \end{equation}
    Solving solution samples $\left\{ {\bf g}\left( t_i, \theta^{\left( q \right)} \right) \right\}_{q=1}^{n_s}$ for all possible sample realizations $\left\{ \theta^{\left( q \right)} \right\}_{q=1}^{n_s}$, we can further calculate random sample realizations of the stochastic solution ${\bm u}\left( {\bf x}, t, \theta \right)$ using \eqref{eq:u_rom}.
    Further, the PC expansion can be used to approximate ${\bm u}\left( {\bf x}, t, \theta \right)$ if an explicit representation is required in practice, which is beyond the scope of this paper and will be investigated in future research.

\subsection{Algorithm flowchart} \label{section:Alg}

\begin{algorithm}[hp]
    \caption{Stochastic LATIN algorithm for solving stochastic elastoplastic problems} \label{Alg_Sto_LATIN}
    \begin{algorithmic}[1]
        
    \State Pre-calculate the deterministic matrices $\left\{ {\bf K}_i \right\}_{i=1}^m$ using \eqref{eq:K_exp}
    \label{alg_LATIN:step_01}

    \While {$\epsilon _{{\bf u}, k} \ge \epsilon _{{\bf u}}$}
    \label{alg_LATIN:step_02}
   
    \State Initialize the random samples $\lambda_k^{\left( 0 \right)}\left( {\widehat{\bm \theta}} \right) \in \mathbb{R}^{n_{s,1}}$ and the time function $g_k^{\left( j \right)}\left( t \right) \in \mathbb{R}^{1 \times n_t}$
    \label{alg_LATIN:step_03}
    
    \While {$\epsilon _{{\bf d}, j} \ge \epsilon _{{\bf d}}$}
    \label{alg_LATIN:step_04}

\vspace{-0.6em}
\textcolor{red}{\rule[4pt]{13cm}{0.05em}}
\vspace{-1.0em}

    \State \hspace{-0.4cm} ${\bm{Global~stage}}:$
    \label{alg_LATIN:step_05}

    \State Update the matrix ${\bf K}_k^{\left( j \right)}$ by \eqref{eq:KdF_ex} and assemble the vector ${\bf F}_k^{\left( j \right)}$ by \eqref{eq:KdF_F}
    \label{alg_LATIN:step_06}
    
    \State Solve the deterministic vector ${\bf d}_k^{\left( j \right)}$ via \eqref{eq:KdF}
    \label{alg_LATIN:step_07}
        
    \State Orthogonalize ${\bf d}_k^{\left( j \right)} \bot \left\{ {\bf d}_i \right\}_{i=1}^{k-1}$ and normalize $\left\| {\bf d}_k^{\left( j \right)} \right\|_2^2 = 1$
    \label{alg_LATIN:step_08}

\vspace{-0.6em}
\textcolor{red}{\rule[4pt]{13cm}{0.05em}}
\vspace{-0.6cm}

\textcolor{blue}{\rule[4pt]{13cm}{0.05em}}
\vspace{-1.0em}

    \State \hspace{-0.4cm} ${\bm{Local~stage}}:$
    \label{alg_LATIN:step_09}
    
    \For {$t_i=t_1, \cdots, t_{n_t}$}
    \label{alg_LATIN:step_10}
    
    \State Initialize the stochastic solution $g_k^{\left( j, 0 \right)}\left( t_i, {\widehat{\bm \theta}} \right) \in \mathbb{R}^{n_{s,1}}$ of time step $t_i$
    \label{alg_LATIN:step_11}
    
    \While {$\epsilon _{g, m}\left( t_i \right) \ge \epsilon _g$}
    \label{alg_LATIN:step_12}
        
    \State Calculate the sample vectors  $a_k^{\left( j, m \right)}\left( {{t_i},{{\widehat{\bm{\theta}}}} } \right) \in \mathbb{R}^{n_{s,1}}$ and $b_k^{\left( j, m \right)}\left( {{t_i},{{\widehat{\bm{\theta}}}} } \right) \in \mathbb{R}^{n_{s,1}}$ by Eqs.~(\ref{eq:g_nl_a}) and (\ref{eq:g_nl_b}) and the elastoplastic constitutive evolution in Section \ref{subsection:SEP_time}
    \label{alg_LATIN:step_13}

    \State Solve the sample vector $\Delta g_k^{\left( j,m \right)}\left( t_i, {\widehat{\bm \theta}} \right) \in \mathbb{R}^{n_{s,1}}$ by \eqref{eq:g_sol_sam}
    \label{alg_LATIN:step_14}

    \State Update the stochastic solution $g_k^{\left( j,m+1 \right)}\left( t_i, {\widehat{\bm \theta}} \right) \in \mathbb{R}^{n_{s,1}}$ by \eqref{eq:g_update}
    \label{alg_LATIN:step_15}
    
    \State Compute the iterative error $\epsilon _{g, m}\left( t_i \right)$, $m \leftarrow m+1$
    \label{alg_LATIN:step_16}

    \EndWhile
    \State {\bf{end}}
    \label{alg_LATIN:step_17}
        
    \EndFor
    \State {\bf{end}}
    \label{alg_LATIN:step_18}
    
    \State Decompose the stochastic solution $g_k^{\left( j \right)}\left( t, {\widehat{\bm \theta}} \right) \in \mathbb{R}^{n_{s,1} \times n_t}$ by \eqref{eq:g_lam_sam}
    \label{alg_LATIN:step_19}

\vspace{-0.6em}
\textcolor{blue}{\rule[4pt]{13cm}{0.05em}}
\vspace{-1.0em}

    \State \hspace{-0.4cm} Compute the iterative error $\epsilon _{{\bf d}, j}$, $j \leftarrow j+1$
    \label{alg_LATIN:step_20}
    
    \EndWhile
    \State {\bf{end}}
    \label{alg_LATIN:step_21}

    \State Update the stochastic solution ${\bf u}_k\left( t,\theta \right) = {\bf u}_{k-1}\left( t,\theta \right) + \lambda_k\left( \theta \right) g_k\left( t \right) {\bf d}_k$ 
    \label{alg_LATIN:step_22}

    \State Compute the iterative error $\epsilon _{{\bf u}, k}$, $k \leftarrow k+1$
    \label{alg_LATIN:step_23}

    \EndWhile
    \State {\bf{end}}
    \label{alg_LATIN:step_24}

    \State ${\bm{Update~stage}}:$ Recalculate the stochastic solution ${\bf g}\left( t, \theta \right)$ by Algorithm \ref{Alg_Newton} for $n_{s,2}$ samples
    \label{alg_LATIN:step_25}
        
    \end{algorithmic}
\end{algorithm}

\begin{algorithm}[!ht]
    \caption{Algorithm of the reduced-order Newton iteration} \label{Alg_Newton}
    \begin{algorithmic}[1]

    \State {\bf Input}: The deterministic matrix ${\bf D} = \left[ {\bf d}_1, \cdots, {\bf d}_k \right] \in \mathbb{R}^{n \times k}$
    \label{alg_NR:step_01}

    \For {$q = 1, \cdots, n_{s,2}$}
    \label{alg_NR:step_02}

    \For {$t_i = t_1, \cdots, t_{n_t}$}
    \label{alg_NR:step_03}
        
        \State Initialize the solution ${\bf{g}}^{\left( 0 \right)}\left( t_i, {{\theta}^{\left( q \right)}} \right) = {\bf{g}}\left( t_{i-1}, {{\theta}^{\left( q \right)}} \right) \in \mathbb{R}^{k}$
        \label{alg_NR:step_04}
        
        \While {$\epsilon _{{\rm NR}, m}\left( {{\theta}^{\left( q \right)}} \right) \ge \epsilon _{{\rm NR}}$}
        \label{alg_NR:step_05}
        
        \State Assemble the reduced-order tangent matrix ${\bf k}_{\rm T}\left( {\bf g}^{\left( m \right)}\left( t_i, \theta^{\left( q \right)} \right), \theta^{\left( q \right)} \right) \in \mathbb{R}^{k \times k}$ and the vector ${\bf f}_{\rm T}\left( {\bf g}^{\left( m \right)}\left( t_i, \theta^{\left( q \right)} \right), \theta^{\left( q \right)} \right) \in \mathbb{R}^{k}$ by Eqs.~(\ref{eq:k_T_rom}) and (\ref{eq:f_rom})
        \label{alg_NR:step_06}
        
        \State Solve the increment $\Delta {\bf g}^{\left( m \right)}\left( t_i, {{\theta}^{\left( q \right)}} \right) \in \mathbb{R}^{k}$ by \eqref{eq:kgf_rom}
        \label{alg_NR:step_07}
        
        \State Update the solution ${\bf{g}}^{\left( m+1 \right)}\left( t_i, {{\theta}^{\left( q \right)}} \right) = {\bf{g}}^{\left( m \right)}\left( t_i, {{\theta}^{\left( q \right)}} \right) + \Delta {\bf{g}}^{\left( m \right)}\left( t_i, {{\theta}^{\left( q \right)}} \right) \in \mathbb{R}^{k}$
        \label{alg_NR:step_08}
        
        \State Calculate the iterative error $\epsilon _{{\rm NR}, m}\left( {{\theta}^{\left( q \right)}} \right)$, $m \leftarrow m+1$
        \label{alg_NR:step_09}
        
        \EndWhile
        \State {\bf{end}}
        \label{alg_NR:step_10}

        \State Update the thermodynamical force ${\bm \beta}\left( {\bf x}, t_i, \theta^{\left( q \right)} \right)$ and the plastic strain ${\bm \varepsilon} _{{\rm p}}\left( {\bf x}, t_i, \theta^{\left( q \right)} \right)$ using Eqs.~(\ref{eq:beta_i}) and (\ref{eq:strin_i})
        \label{alg_NR:step_11}
        
        \EndFor
        \State {\bf{end}}
        \label{alg_NR:step_12}

        \EndFor
        \State {\bf{end}}
        \label{alg_NR:step_13}
        
    \end{algorithmic}
\end{algorithm}

    The proposed method is summarized in Algorithm \ref{Alg_Sto_LATIN}.
    The deterministic matrices $\left\{ {\bf K}_i \right\}_{i=1}^m$ are assembled by \eqref{eq:K_exp} in step \ref{alg_LATIN:step_01}, which is only performed once for all possible global stages.
    It is noted that an \eqref{eq:C_expansion}-like affine approximation of random inputs may need to be pre-performed in practice.
    Two loops are then used to solve the stochastic solution, including the outer loop from step \ref{alg_LATIN:step_02} to \ref{alg_LATIN:step_24} and the inner loop from step \ref{alg_LATIN:step_04} to \ref{alg_LATIN:step_21}.
    To calculate each triplet $\left\{ \lambda_k\left( \theta \right), g_k\left( t \right), {\bm d}_k\left( {\bf x} \right) \right\}$ using the inner loop, a random sample vector $\lambda_k^{\left( 0 \right)}\left( {\widehat{\bm \theta}} \right) \in \mathbb{R}^{n_{s,1}}$ and a discretized time vector $g_k^{\left( j \right)}\left( t \right) \in \mathbb{R}^{1 \times n_t}$ are initialized in step \ref{alg_LATIN:step_03}.
    Any nonzero vectors of size $n_{s,1}$ can be used as the initialization for $\lambda_k\left( \theta \right)$ and the sample size $n_{s,1}$ can be a small number during the iteration process, which has little influence on the accuracy of the final stochastic solution, but saves a lot of computational effort. 
    However, the choice of $n_{s,1}$ is currently trial and more detailed research is required.
    
    The stochastic LATIN iteration is then performed in the inner loop.
    The global stage is from step \ref{alg_LATIN:step_05} to \ref{alg_LATIN:step_08} and used to solve the spatial vector ${\bf d}_k^{\left( j \right)} \in \mathbb{R}^n$.
    For each global stage, we only need to update the coefficients $\left\{ z_{k,i}^{\left( j \right)} \right\}_{i=1}^m$ by \eqref{eq:KdF_ex}, which is very cheap due to the reuse of deterministic matrices $\left\{ {\bf K}_i \right\}_{i=1}^m$.
    The computational effort is thus mainly concentrated on the assembly of the vector ${\bf F}_k^{\left( j \right)} \in \mathbb{R}^n$ by \eqref{eq:KdF_F}, whose cost can be reduced to some extent by adopting a small sample size $n_{s,1}$.
    To speed up the convergence and achieve a small retained number $k$, we let the vector ${\bf d}_k^{\left( j \right)}$ orthogonal to the obtained vector $\left\{ {\bf d}_i \right\}_{i=1}^{k-1}$ during the whole iterative process, which is achieved by the Gram-Schmidt orthonormalization

    \vspace{-1.0em}
    \begin{equation} \label{eq:GS_Or}
        {\bf d}_k^{\left( j \right)} = {\bf d}_k^{\left( j \right)} - \sum\limits_{i=1}^{k-1} \frac{ {\bf d}_k^{\left( j \right){\rm T}}{\bf d}_i }{{\bf d}_i^{\rm T}{\bf d}_i} {\bf d}_i, \quad {\bf d}_k^{\left( j \right){\rm T}} {\bf d}_k^{\left( j \right)} = 1.
    \end{equation}
    
    Based on the known spatial vector ${\bf d}_k^{\left( j \right)}$, the local stage is then performed from step \ref{alg_LATIN:step_09} to \ref{alg_LATIN:step_19} and used to calculate the temporal-stochastic coupling solution $g_k^{\left( j \right)}\left( t, {\widehat{\bm \theta}} \right) \in \mathbb{R}^{n_{s,1} \times n_t}$.
    For each time step $t_i$, a stochastic Newton iteration from step \ref{alg_LATIN:step_12} to \ref{alg_LATIN:step_17}
    is adopted to solve local nonlinear stochastic problems.
    The sample vectors $a_k^{\left( j, m \right)}\left( {{t_i},{{\widehat{\bm{\theta}}}} } \right) \in \mathbb{R}^{n_{s,1}}$ and $b_k^{\left( j, m \right)}\left( {{t_i},{{\widehat{\bm{\theta}}}} } \right) \in \mathbb{R}^{n_{s,1}}$ are obtained in step \ref{alg_LATIN:step_13} by assembling each random sample realization, which also involves the temporal evolution of the stochastic elastoplastic constitutive model.
    Following the iterations of all time steps, we perform rank-1 SVD for the temporal-stochastic solution $g_k^{\left( j \right)}\left( t, {\widehat{\bm \theta}} \right) \in \mathbb{R}^{n_{s,1} \times n_t}$ in step \ref{alg_LATIN:step_19}.
	After the inner loop, the outer loop is to recursively update the stochastic solution ${\bf u}_k\left( t,\theta \right)$ in step \ref{alg_LATIN:step_22} with the new triplet $\left\{ \lambda_k\left( \theta \right), g_k\left( t \right), {\bf d}_k \right\}$.
	Based on the spatial reduced basis functions, a recalculation process of the temporal-stochastic coupling solution vector ${\bf g}\left( t, \theta \right) \in \mathbb{R}^k$ is executed in step \ref{alg_LATIN:step_25} by using Algorithm \ref{Alg_Newton}.
    It is noted that the sample size $n_{s,1}$ used for the iterative process is much smaller than the sample size $n_{s,2}$ used for the recalculation process, which reduces much computational effort of the iterative process.

    Four error indicators are involved in the above algorithms.
    For each time step of the Newton iteration in the local stage, the error indicator $\epsilon _{g, m}\left( t_i \right)$ in step \ref{alg_LATIN:step_16} is given by

    \vspace{-1.0em}
    \begin{equation} \label{eq:It_err_g}
        \epsilon _{g,m}\left( t_i \right) = \frac{\left\| \Delta g_k^{\left( j,m \right)}\left( t_i, {\widehat{\bm \theta}} \right) \right\|_2^2}{\left\| g_k^{\left( j,m \right)}\left( t_i, {\widehat{\bm \theta}} \right) \right\|_2^2} ,
    \end{equation}
    where $\left\| \cdot \right\|_2$ denotes the $L_2$ norm. 
    This error measures the contribution of the $m$-th increment $\Delta g_k^{\left( j,m \right)}\left( t_i, {\widehat{\bm \theta}} \right)$ to the solution $g_k^{\left( j,m \right)}\left( t_i, {\widehat{\bm \theta}} \right)$.
    Similarly, the error indicator $\epsilon_{{\rm NR},m}\left( {{\theta}^{\left( q \right)}} \right)$ of the reduced-order Newton iteration in step \ref{alg_NR:step_09} of Algorithm \ref{Alg_Newton} is given by

	\vspace{-1.0em}
	\begin{equation}\label{eq:err_NR}
	   \epsilon_{{\rm NR},m}\left( {{\theta}^{\left( q \right)}} \right) = \frac{ \left\| \Delta {\bf{g}}^{\left( m \right)}\left( t, {{\theta}^{\left( q \right)}} \right) \right\|_2^2 }{\left\| {\bf{g}}^{\left( m \right)}\left( t, {{\theta}^{\left( q \right)}} \right) \right\|_2^2}.
	\end{equation}
    The error indicator $\epsilon _{{\bf d},j}$ of the inner loop in step \ref{alg_LATIN:step_20} of Algorithm \ref{Alg_Sto_LATIN}, i.e., the stopping criterion for each triplet, is given by

    \vspace{-1.0em}
    \begin{equation} \label{eq:It_err_d}
    	\epsilon _{{\bf d},j} = \frac{\left\| {{\bf{d}}^{\left( j \right)}_k} - {{\bf{d}}^{\left( j-1 \right)}_k} \right\|_2^2}{\left\| {{\bf{d}}^{\left( j \right)}_k} \right\|_2^2} = 2 - 2{{\bf{d}}^{\left( j \right) {\rm T}}_k}{{\bf{d}}^{\left( j-1 \right)}_k}.
    \end{equation}

    Further, The iterative error of the outer loop in step \ref{alg_LATIN:step_23}, i.e., the stopping criterion for the stochastic solution ${\bm u}\left( {\bf x}, t, \theta \right)$, is given by

	\vspace{-1.0em}
	\begin{equation} \label{eq:It_err_u}
    	\epsilon _{{\bf u},k} = \frac{{\left\| {{\bf{u}}_{k}}\left( {t,{\widehat{\bm{\theta}}} } \right) - {{\bf{u}}_{k-1}}\left( {t,{\widehat{\bm{\theta}}} } \right) \right\|_2^2}}{{\left\| {{\bf{u}}_{k}}\left( {t,{\widehat{\bm{\theta}}} } \right) \right\|_2^2}} 
    	= \frac{{{\widehat{\mathbb{E}}}\left\{ {\lambda _k^2\left( {{\widehat{\bm{\theta}}} } \right)} \right\}{\bf{d}}_k^{\rm T}{{\bf{d}}_k} \int_0^T g_k^2\left( t \right){\rm d}t}}{{{\sum\limits_{i,j = 1}^k {{\widehat{\mathbb{E}}}\left\{ {{\lambda _i}\left( {{\widehat{\bm{\theta}}} } \right){\lambda _j}\left( {{\widehat{\bm{\theta}}} } \right)} \right\}{\bf{d}}_i^{\rm T}{{\bf{d}}_j} \int_0^T g_i\left( t \right)g_j\left( t \right){\rm d}t } } } }
    	= \frac{{{\widehat{\mathbb{E}}}\left\{ {\lambda _k^2\left( {{\widehat{\bm{\theta}}} } \right)} \right\}}}{{\sum\limits_{i = 1}^k {{\widehat{\mathbb{E}}}\left\{ {\lambda _i^2\left( {{\widehat{\bm{\theta}}} } \right)} \right\}} }},
	\end{equation}
	which measures the contribution of the $k$-th triplet to the stochastic solution ${\bm u}_k\left( {\bf x}, t, \theta \right)$.
    However, since the random variables $\left\{ \lambda_i\left( \theta \right) \right\}_{i=1}^k$ are solved in a sequential way, this error definition does not keep decreasing as the retained term $k$ increases. 
    To avoid this issue, we replace $\left\{ \lambda_i\left( \theta \right) \right\}_{i=1}^k$ with a set of equivalent random variables $\left\{ {\overline \lambda}_i\left( \theta \right) \right\}_{i=1}^k$  \cite{zheng2022weak, zheng2023stochastic}.
    To this end, let us consider the autocorrelation matrix of the random vector ${\bf \Lambda}\left( \theta \right) = \left[ \lambda_1\left( \theta \right), \cdots, \lambda_k\left( \theta \right) \right]^{\rm T} \in \mathbb{R}^k$ given by ${\bf Cov}_{\bf{\Lambda}} = \mathbb{E} \left\{ {\bf{\Lambda}}\left( \theta \right) {\bf{\Lambda}}\left( \theta \right)^{\rm T} \right\}$, which can be decomposed into
	
	\vspace{-1.0em}
	\begin{equation}\label{eq:dom_cov}
	   {\bf Cov}_{\bf{\Lambda}} = {\bf{QZQ}}^{\rm T}
	\end{equation}
	by the eigendecomposition, where ${\bf Q} \in \mathbb{R}^{k \times k}$ is an orthonormal matrix and ${\bf Z} \in \mathbb{R}^{k \times k}$ is a diagonal matrix consisting of descending eigenvalues of ${\bf Cov}_{\bf{\Lambda}}$. 
    We can further rewrite the stochastic solution ${\bm u}\left( {\bf x}, t, \theta \right)$ in \eqref{eq:u_app} as
	
	\vspace{-1.0em}
	\begin{equation}\label{eq:dom_u}
	   {\bm u}\left( {\bf x}, t, \theta \right) = {\bm D}\left( {\bf x} \right) {\bm G}\left( t \right) \underbrace{{\bf{QQ}}^{\rm T}}_{{\bf I}_k} {\bf{\Lambda }}\left( \theta  \right),
	\end{equation}
    where the diagonal matrix ${\bm G}\left( t \right) = {\rm diag}\left[ g_1\left( t \right), \cdots, g_k\left( t \right) \right] \in \mathbb{R}^{k \times k}$.
	We let an equivalent random vector be ${\bf{\overline \Lambda }}\left( \theta  \right) = {\bf{Q}}^{\rm T} {\bf{\Lambda }}\left( \theta  \right) = \left[ {{\overline \lambda }_1}\left( \theta  \right), \cdots,{{\overline \lambda }_k}\left( \theta  \right) \right]^{\rm T} \in \mathbb{R}^k$,
	and its autocorrelation function is ${{\bf Cov}}_{\bf{\overline \Lambda}} = {\bf Z}$.
    Therefore, \eqref{eq:It_err_u} is reformulated as
	
	\vspace{-1.0em}
	\begin{equation}\label{eq:err_glo}
    	\epsilon _{{\bf u},k} 
    	= \frac{{\mathbb{E}\left\{ {{\overline \lambda} _k^2\left( \theta  \right)} \right\}}}{{\sum\limits_{i = 1}^k {\mathbb{E}\left\{ {{\overline \lambda} _i^2\left( \theta  \right)} \right\}} }} = \frac{{\bf Z}_k}{{\rm{Tr}}\left({\bf{Z}} \right)},
	\end{equation}
	where ${\rm{Tr}}\left( \cdot \right)$ is the trace operator and ${\bf Z}_k$ is the $k$-th diagonal element of the matrix ${\bf Z}$. 
	In this way, the iterative error $\epsilon_{{\bf u},k}$ keeps decreasing as the retained item $k$ increases. It is noted that \eqref{eq:dom_u} is equivalent to \eqref{eq:u_app} and only provides a new representation of the random vector ${\bf{\Lambda }}\left( \theta  \right)$. 

\section{Numerical examples} \label{section:Examples}
    In this section, four numerical examples are presented to illustrate performance of the proposed stochastic LATIN method.
    For all cases, the convergence criteria for inner and outer loops of Algorithm \ref{Alg_Sto_LATIN} are set as $\epsilon _{\bf d} = 1 \times 10^{-3}$ and $\epsilon _{\bf u} = 1 \times 10^{-8}$, respectively.
    The convergence criterion for Newton iterations in both the local stage of Algorithm \ref{Alg_Sto_LATIN} and Algorithm \ref{Alg_Newton} is set as $\epsilon _g = \epsilon _{\rm NR} = 1 \times 10^{-10}$.
    The sample size $n_{s,2}$ during the iterative process of Algorithm \ref{Alg_Sto_LATIN} is 10 times the dimension of random and/or parameter inputs.
    Further, reference solutions are obtained by $1 \times 10^3$ Newton-based standard Monte Carlo simulations.
    To eliminate the influence from sampling processes, the same random sample realizations are also applied to the update stage of the proposed stochastic LATIN method.
    All numerical tests are performed on one core of a desktop computer (24 cores, Intel Xeon Silver 4116 CPU, 2.10GHz).

\subsection{Problem setting}

\begin{figure}[!h]
    \begin{minipage}{0.587\textwidth}
        \centering
        \subfloat[][Physical model and its finite element mesh.]{\label{fig_Model} \includegraphics[width=0.9\textwidth]{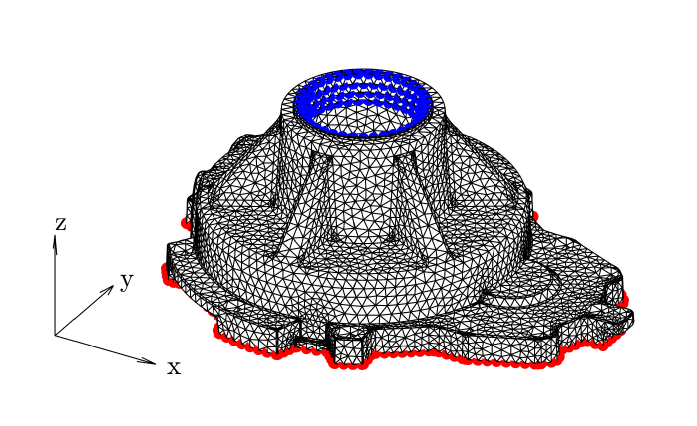}}
    \end{minipage}
    \begin{minipage}{0.413\textwidth}
        \centering
        \subfloat[][Top view and the reference point $A$.]{\label{fig_Model_top} \includegraphics[width=0.9\textwidth]{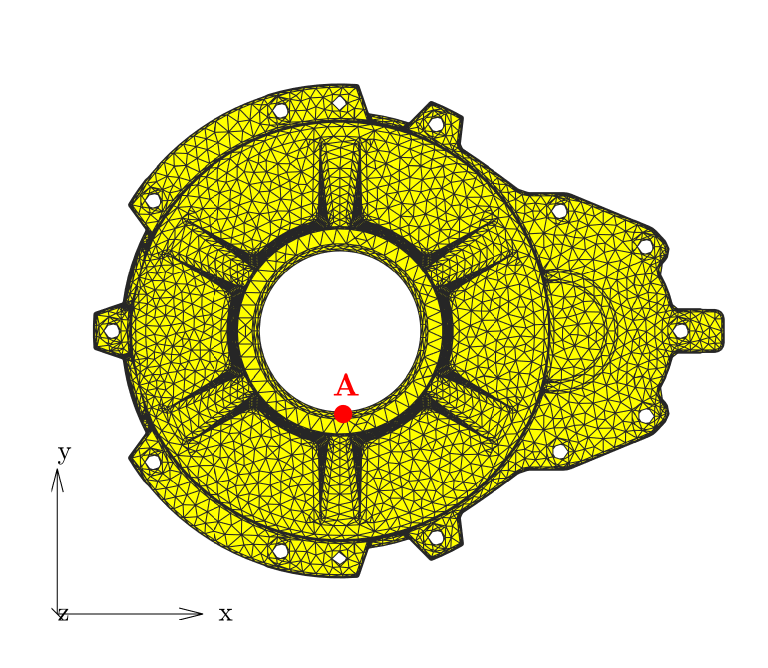}}
    \end{minipage}
    \begin{minipage}{1.0\textwidth}
        \centering
        \subfloat[][Fixed boundaries.]{\label{fig_Model_bottom} \includegraphics[width=0.495\textwidth]{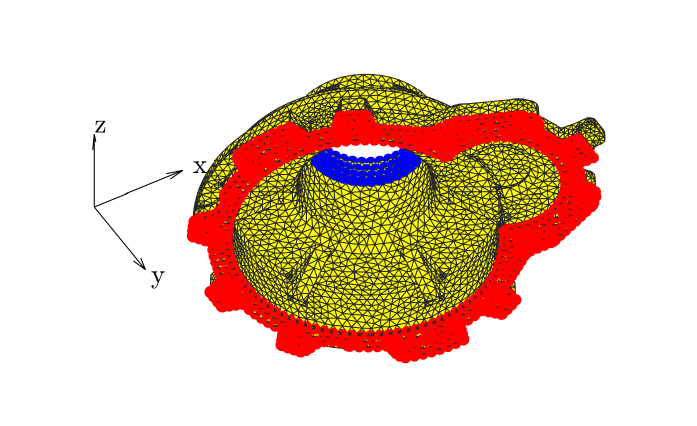}}
    \end{minipage}
    \caption{Gearbox cover model and its finite element mesh.} \label{fig_Model_Mesh}
\end{figure}

\begin{figure}[ht]
        \centering
        \includegraphics[width=0.5\linewidth]{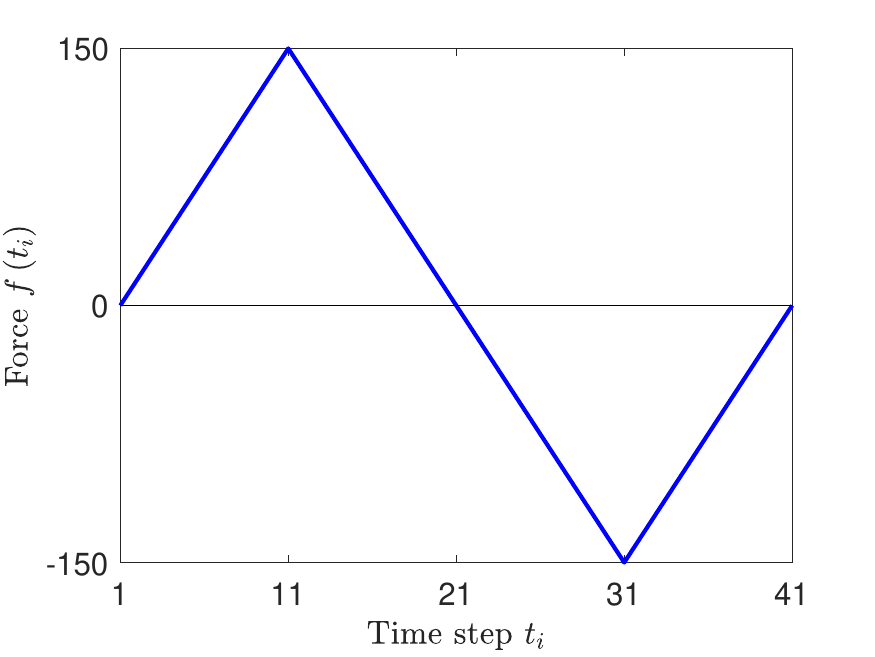}
        \caption{The history of cyclic load.}
        \label{fig_Force}
    \end{figure}

    We consider the stochastic and/or parameterized elastoplastic analysis for a 3D gearbox cover shown in \figref{fig_Model}.
    A cyclic force $f\left( x, y, z, t \right)$, as shown in \figref{fig_Force}, is applied distributedly in the blue zone (seen from \figref{fig_Model}) along the $y$ direction, and its magnitude is $150$~${\rm N/mm^2}$.
    The time variable is discretized as 41 time steps $t \in \left[ t_1, \cdots, t_{41} \right]$ for all examples under consideration.
    The model is fixed on the red boundary (shown in \figref{fig_Model_bottom}), that is, the corresponding three displacement components are $u_x\left( x, y, z, \theta \right) = u_y\left( x, y, z, \theta \right) = u_z\left( x, y, z, \theta \right) = 0$.
    The finite element mesh is discretized using the linear tetrahedron element, with a total of 10072 nodes, 32284 elements and 30216 spatial degrees of freedom.

    The hardening parameter in \eqref{eq:inter_pat} is considered as a deterministic value $\varpi = 1 \times 10^4$~Pa, but it is easy to be extended to a parameter or a random variable/field.
	The stochastic material properties $\kappa\left( x, y, z, \theta \right)$ and $\mu\left( x, y, z, \theta \right)$ in \eqref{eq:Ce} are given by
		
	\vspace{-1.0em}
	\begin{equation} \label{eq:Lame}
    	\kappa\left( x, y, z, \theta \right) = \frac{E\left( x, y, z, \theta \right)}{3\left( 1 - 2\nu \right)}, \quad
    	\mu\left( x, y, z, \theta \right) = \frac{E\left( x, y, z, \theta \right)}{2\left( 1 + \nu \right)},
	\end{equation}
	where the Poisson rate $\nu = 0.29$ and the Young's modulus $E\left( x,y,\theta \right)$ will be respectively modeled as a parameter, a random variable and a random field in subsequent sections.
    Also, the yield stress $\sigma_{\rm Y}\left( x, y, z, \theta \right )$ will be modeled similarly.

\subsection{Example 1: Stochastic elastoplastic problem} \label{subsection:Sto_Case}
    In this example, we consider a stochastic elastoplastic problem.
    Specifically, the Young's modulus and the yield stress are modeled as two Gaussian random variables

    \vspace{-1.0em}
    \begin{equation*}
        E\left( \theta \right) \sim {\cal N}\left( m_{E}, \chi_E^2 \right), \quad 
        \sigma_{\rm Y}\left( \theta \right) \sim {\cal N}\left( m_{\sigma_{\rm Y}}, \chi_{\sigma_{\rm Y}}^2 \right),
    \end{equation*}
    where the mean values are $m_E = 2.11 \times 10^5$~MPa, $m_{\sigma_{\rm Y}} = 245$~MPa and the standard deviations are $\chi_E = 0.1 m_E$, $\chi_{\sigma_{\rm Y}} = 0.1 m_{\sigma_{\rm Y}}$.
    It is noted that to ensure the physical fact that material properties are positive, the random events $\theta^*$ such that $E\left( \theta^* \right) < 0.1 m_E$ and $\sigma_{\rm Y}\left( \theta^* \right) < 0.1 m_{\sigma_{\rm Y}}$ are dropped out in practical numerical implementations.

\begin{figure}[ht]
    \centering
    \includegraphics[width=0.5\linewidth]{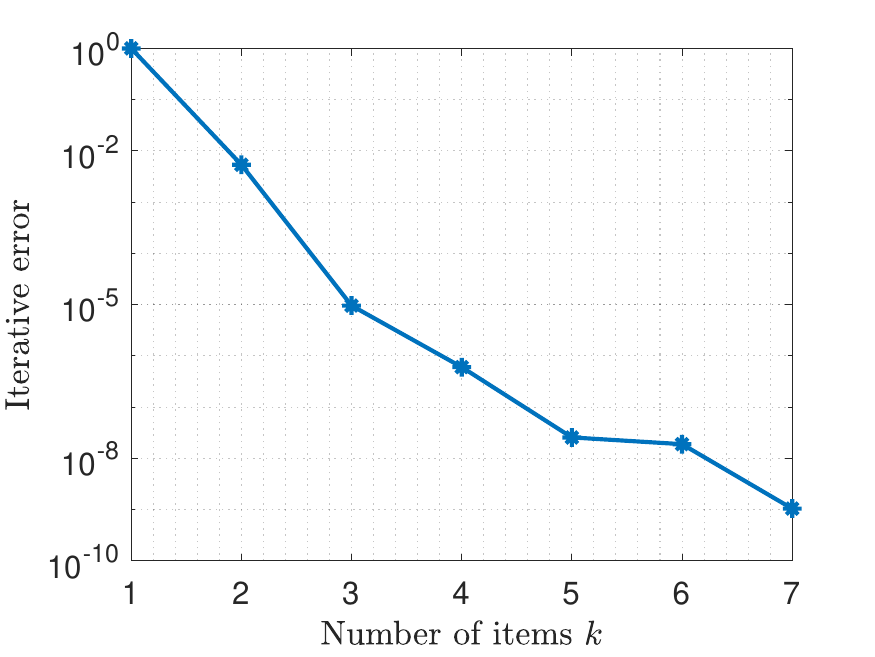}
    \caption{Iterative errors of diﬀerent retained terms.}
    \label{fig_Ex_1_It_err}
\end{figure}

    For numerical implementations of the proposed stochastic LATIN method, the sample size $n_{s,1} = 20$ is adopted for the iterative process of Algorithm \ref{Alg_Sto_LATIN}.
    In our experience, these sample realizations of the random inputs can be generated without careful design and have little influence on the proposed method.
    They only need to roughly follow the given probability distributions.
    Iterative errors of different diﬀerent retained terms $k$ are shown in \figref{fig_Ex_1_It_err}.
    It is seen that only 7 retained terms are required to achieve the specified precision, which demonstrates the good convergence of the proposed stochastic LATIN method.
    The total number of linear global equations (\ref{eq:KdF}) needed to be solved during the iterative process is only 26.
    If the Newton iteration-based method is used, 41 nonlinear stochastic problems need to be solved, which is even more than the number of linear deterministic equations from the proposed method.
    Therefore, the proposed method is much more efficient than the Newton iteration-based methods. 

\begin{figure}[ht]
    \begin{minipage}{0.5\textwidth}
        \centering
        \subfloat[][The stochastic displacements obtained by the MCS and the stochastic LATIN method without and with update.]{\label{fig_e1_Dis_A_01} \includegraphics[width=1.0\textwidth]{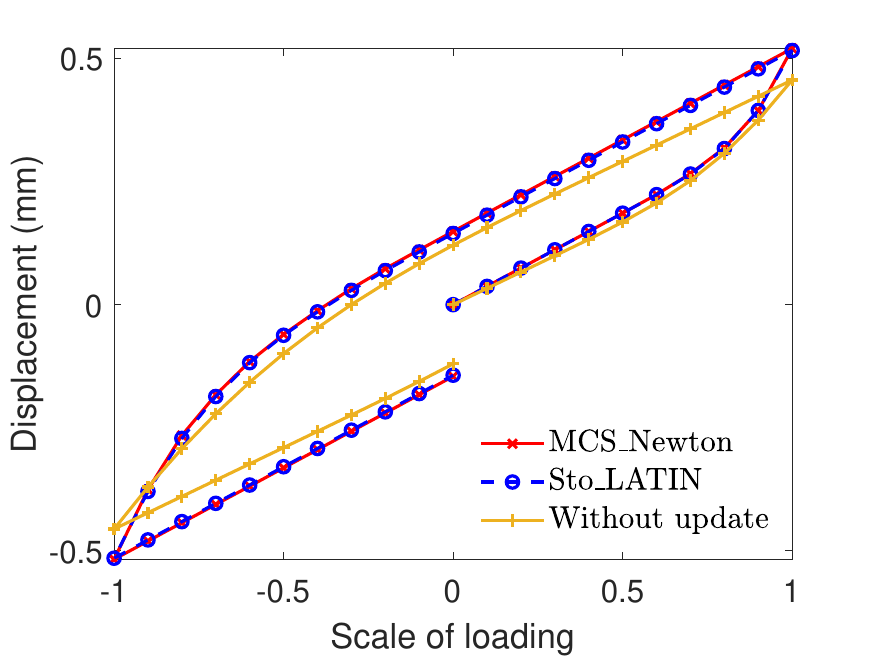}}
    \end{minipage}
    \hspace{0.12cm}
    \begin{minipage}{0.5\textwidth}
        \centering
        \subfloat[][Absolute errors of the stochastic displacements referring to the Newton-based MCS reference solution.]{\label{fig_e1_Dis_A_01_Ab_Err} \includegraphics[width=1.0\textwidth]{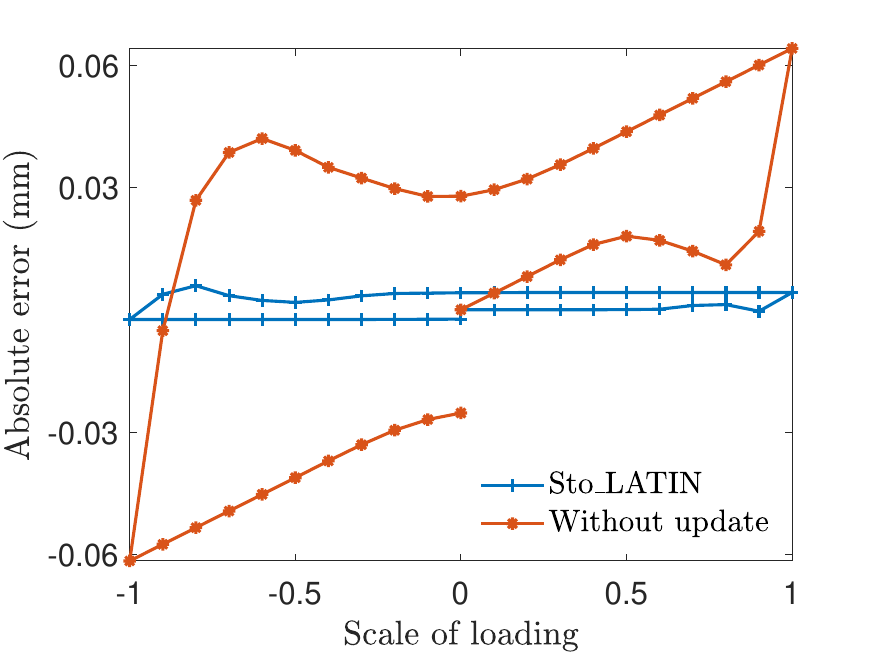}}
    \end{minipage}
    \caption{The stochastic displacements of the reference point $A$ in the $y$ direction of the sample realization $\left\{ E\left( \theta^* \right), \sigma_{\rm Y}\left( \theta^* \right) \right\} = \left\{ 2.25 \times 10^5, 270.07 \right\}$ obtained by the Newton-based MCS and the stochastic LATIN iterations without and with update (left) and the absolute errors referring to the Newton-based MCS displacement (right).} \label{fig_e1_Dis_A_01_err}
\end{figure}

\begin{figure}[!ht]
    \centering
    \includegraphics[width=0.5\linewidth]{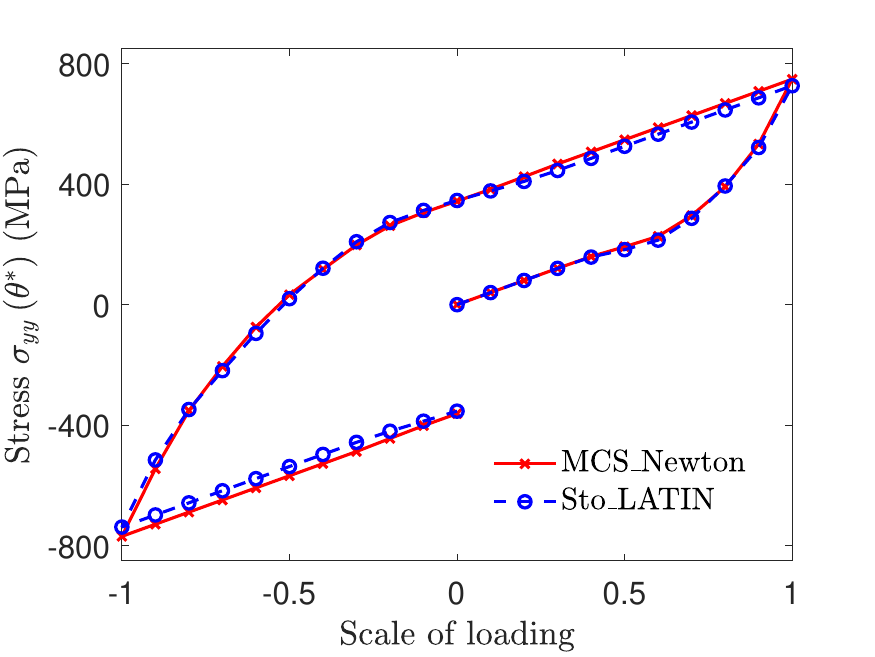}
    \caption{The stochastic stress $\sigma_{yy}\left( \theta^* \right)$ of the point $A$ in the $y$ direction obtained by the MCS and the proposed stochastic LATIN method.}
    \label{fig_e1_stress}
\end{figure}

    To test the computational accuracy of the proposed method, we consider the stochastic solution concerning the sample realization $\left\{ E\left( \theta^* \right), \sigma_{\rm Y}\left( \theta^* \right) \right\} = \left\{ 2.25 \times 10^5, 270.07 \right\}$, which is also one of the 20 sample realizations used for the iterative process.
    The stochastic displacements in the $y$ direction of the reference point $A$ (shown in \figref{fig_Model_top}) are shown in \figref{fig_e1_Dis_A_01}, where \textit{MCS\_Newton}, \textit{Sto\_LATIN} and \textit{Without update} represent the solutions obtained by the Newton iteration-based MCS, and the proposed stochastic LATIN method with the update stage in step \ref{alg_LATIN:step_25} of Algorithm \ref{Alg_Sto_LATIN} and without the update stage.
    Their absolute errors referring to the MCS reference solution are plotted in \figref{fig_e1_Dis_A_01_Ab_Err}.
    It is seen that the solution without update is less accurate and the update stage can greatly improve the accuracy.  
    Further, we adopt the $L_2$ norm error to quantify the approximation error, and only the displacement in the $y$ direction is considered here.
    For a given spatial position ${\bf x}^*$ (i.e. the point $A$ here) and a sample realization $\theta^*$, the $L_2$ approximation error is defined over the temporal domain
    
    \vspace{-1.0em}
    \begin{equation} \label{eq:L2_err}
        \epsilon\left( {\bf x}^*, \theta^* \right) = \frac{ \left\| u_y^{\rm MC}\left( {\bf x}^*, t, \theta^* \right) - u_y^{\rm SL}\left( {\bf x}^*, t, \theta^* \right) \right\|_2}{ \left\| u_y^{\rm MC}\left( {\bf x}^*, t, \theta^* \right) \right\|_2 } ,
    \end{equation}
    where $u_y^{\rm MC}\left( {\bf x}^*, t, \theta^* \right)$ and $u_y^{\rm SL}\left( {\bf x}^*, t, \theta^* \right)$ are the stochastic displacements obtained by the MCS and the proposed method, respectively.
    In this way, the $L_2$ errors of the stochastic displacements without and with update are $1.59 \times 10^{-2}$ and $1.10 \times 10^{-4}$, respectively.
    The approximation accuracy is improved by two orders of magnitude by the update stage.
    Furthermore, stochastic stress is typically the quantity of interest in engineering practice.
    Although only stochastic displacement is considered here, it is easy to calculate stochastic stresses from stochastic displacements and the elastoplastic constitutive model.
    The stochastic stresses corresponding to the stochastic displacements in \figref{fig_e1_Dis_A_01} are compared in \figref{fig_e1_stress}.
    The two stochastic stresses are still in good agreement.

\begin{figure}[ht]
    \centering
    \includegraphics[width=1.0\linewidth]{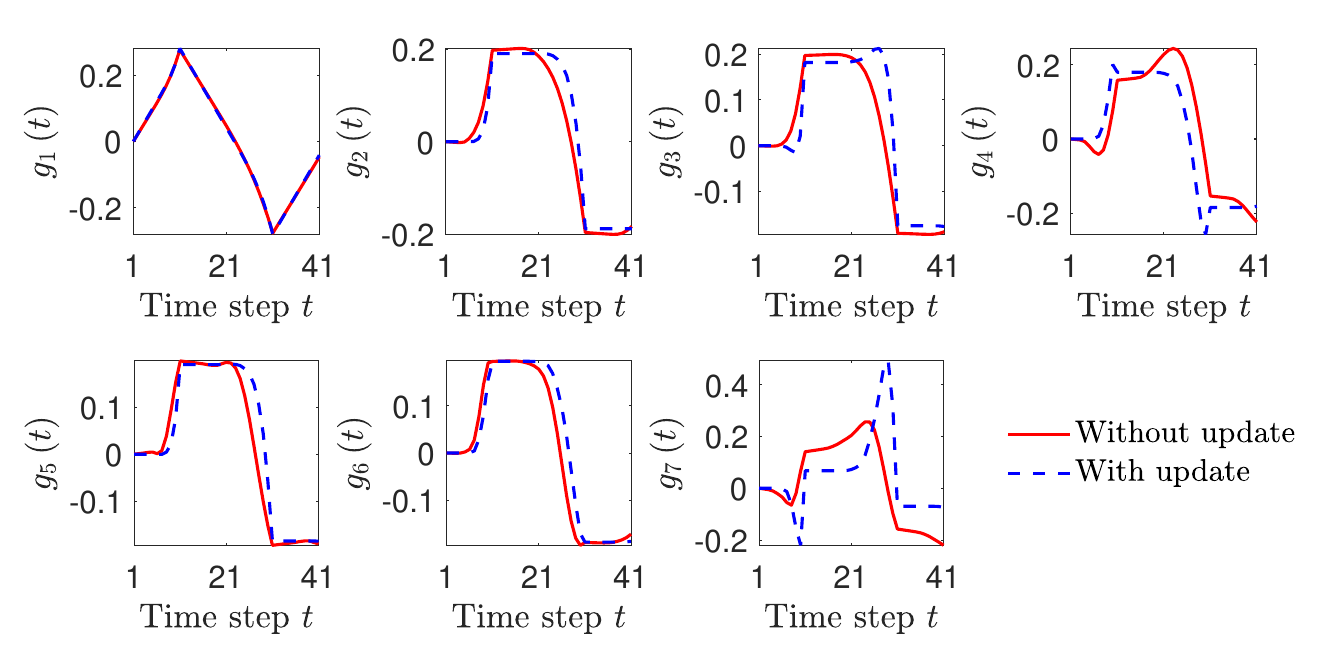}
    \caption{The temporal functions $\left\{ g_i\left( t \right) \right\}_{i=1}^7$ without and with update.}
    \label{fig_e1_Gt}
\end{figure}

\begin{table}[ht]
    \centering
    \caption{Unupdated random variables $\left\{ \lambda_i\left( \theta^* \right) \right\}_{i=1}^7$ and updated random variables $\left\{ \lambda_i^{\rm up}\left( \theta^* \right) \right\}_{i=1}^7$.}\label{tab_e1_lam}
    \begin{tabular}{cccccccc}
        \toprule
         $i$ & 1 & 2 & 3 & 4 & 5 & 6 & 7 \\
        \midrule
        $\lambda_i\left( \theta^* \right)$  & 37.23 & 3.45 & 0.69 & 0.45 & 0.34 & 0.23 & 0.02 \\
            $\lambda_i^{\rm up}\left( \theta^* \right)$  & 39.82 & 5.38 & 1.12 & 0.42 & 0.57 & 0.28 & 0.06 \\
        \bottomrule
    \end{tabular}
\end{table}

\begin{figure}[!ht]
    \begin{minipage}{0.5\textwidth}
        \centering
        \subfloat[][Absolute errors different retained items $k=4 \sim 8$.]{\label{fig_e1_Ab_err_k} \includegraphics[width=1.0\textwidth]{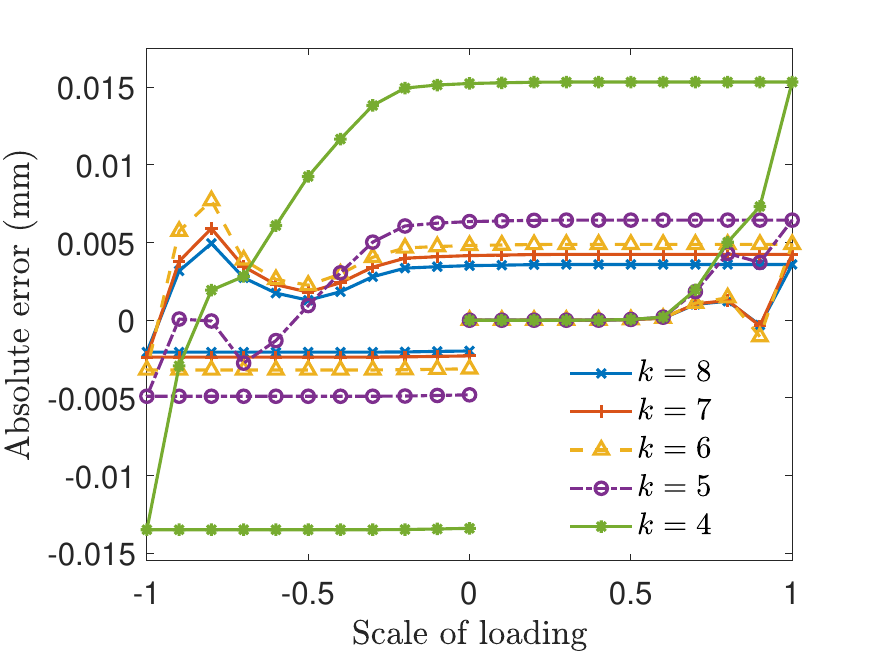}}
    \end{minipage}
    \hspace{0.12cm}
    \begin{minipage}{0.5\textwidth}
        \centering
        \subfloat[][$L_2$ errors different retained items $k=4 \sim 8$.]{\label{fig_e1_L2_err_k} \includegraphics[width=1.0\textwidth]{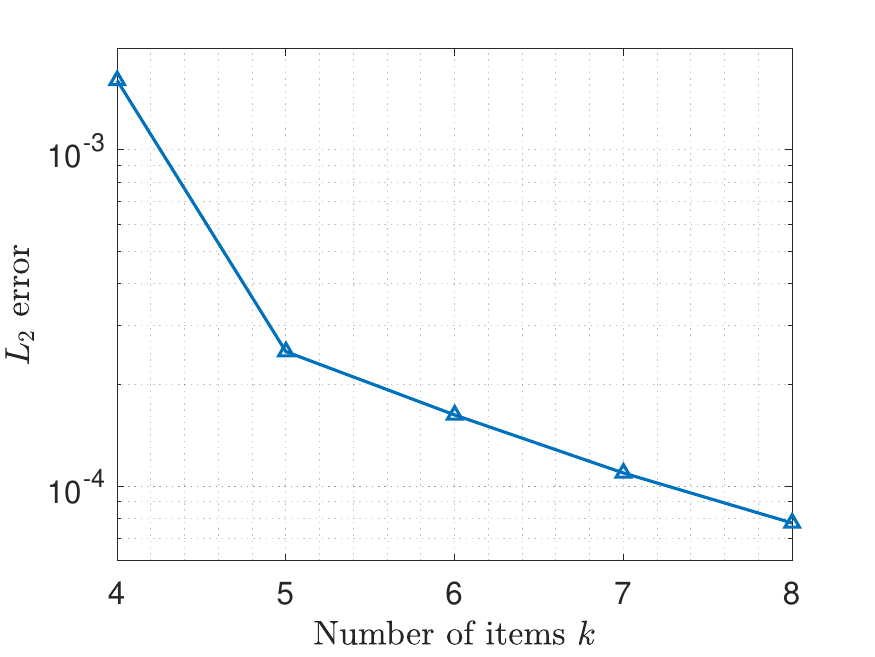}}
    \end{minipage}
    \caption{Absolute errors (left) of the stochastic displacements the point $A$ in the $y$ direction obtained by the stochastic LATIN method with different retained items $k=4 \sim 8$ referring to the MCS reference solution and corresponding $L_2$ errors (right).} \label{fig_e1_Dis_A_sam_Ab_L2_err}
\end{figure}

    The comparison of temporal functions $\left\{ g_i\left( t \right) \right\}_{i=1}^7$ of the above sample input $\left\{ E\left( \theta^* \right), \sigma_{\rm Y}\left( \theta^* \right) \right\} = \left\{ 2.25 \times 10^5, 270.07 \right\}$ without and with update is found in \figref{fig_e1_Gt}, and sample realizations of the corresponding random variables are listed in Table \ref{tab_e1_lam}.
    The temporal functions without update are fixed for all possible sample realizations of $\left\{ \lambda_i\left( \theta \right) \right\}_{i=1}^7$, that is, each sample realization of the stochastic solution is achieved by the combination of $\left\{ g_i\left( t \right) \right\}_{i=1}^7$ and a sample realization $\left\{ \lambda_i\left( \theta^* \right) \right\}_{i=1}^7$.
    The update stage adaptively recalculates the updated temporal functions and the updated sample realizations of the random variables $\left\{ \lambda_i^{\rm up}\left( \theta \right) \right\}_{i=1}^7$, which thus has better approximation accuracy.
    For the sample input under consideration, it is seen that most temporal functions without and with update have similar modes, but the random variables are updated to accommodate the stochastic solution of the given sample input. 
    Further, the approximation accuracy of different retained terms $k = 4,5,6,7,8$ are compared in \figref{fig_e1_Ab_err_k} and corresponding $L_2$ errors are shown in \figref{fig_e1_L2_err_k}.
    We also investigate the errors of the 8-term approximation although the 7-term approximation has achieved the given convergence error.
    It is seen that the approximation accuracy increases as the number of retained terms increases.
    In practice, a higher-accuracy stochastic solution is easily reached by retaining more triplets.

\begin{figure}[ht]
    \centering
    \includegraphics[width=0.5\linewidth]{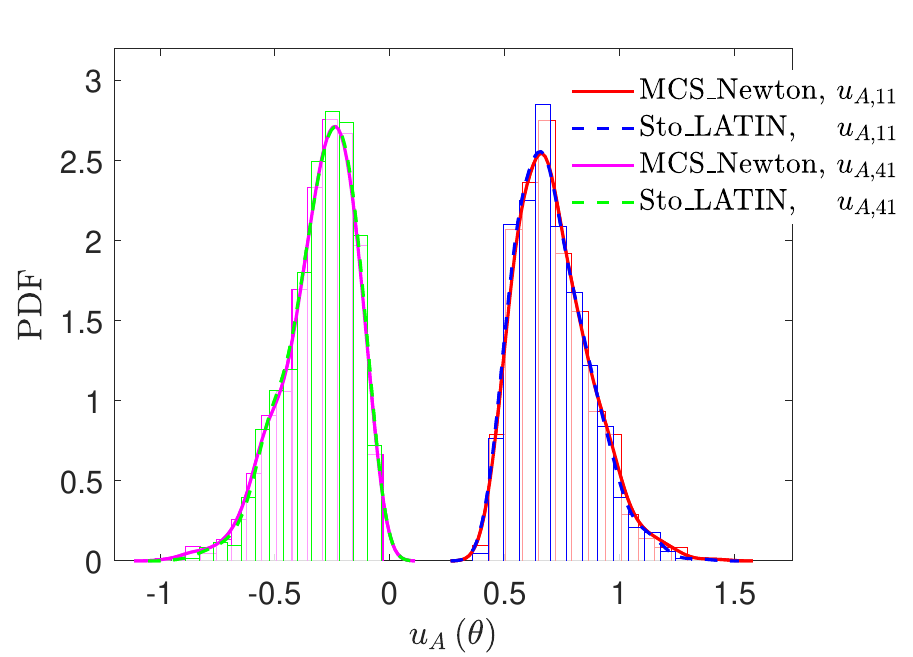}
    \caption{The PDFs of the stochastic displacement of the point $A$ in the $y$ direction at time steps $t_{11}$ and $t_{41}$ obtained by the MCS and the proposed stochastic LATIN method.}
    \label{fig_e1_PDF_A_11_41}
\end{figure}

    The probability density functions (PDFs) of stochastic displacement of the point $A$ in the $y$ direction at the time steps 11 and 41 are compared in \figref{fig_e1_PDF_A_11_41} and they are solved using the MCS and the proposed method with $1 \times 10^3$ random sample realizations.
    For both time steps, the two approaches have good agreements, which demonstrates the high accuracy of the proposed method.
    Due to the expensive computational burden of MCS, only $1 \times 10^3$ sample realizations are executed, which may not be enough to reach well-converged PDFs.
    However, the obtained PDFs are sufficient to illustrate that the proposed stochastic LATIN method has comparable accuracy with the Newton iteration-based MCS.
    
\begin{figure}[ht]
    \begin{minipage}{0.5\textwidth}
        \centering
        \subfloat[][The stochastic displacements of two sample realizations obtained by the MCS and the stochastic LATIN iteration.]{\label{fig_e1_Dis_A_2_sam} \includegraphics[width=1.0\textwidth]{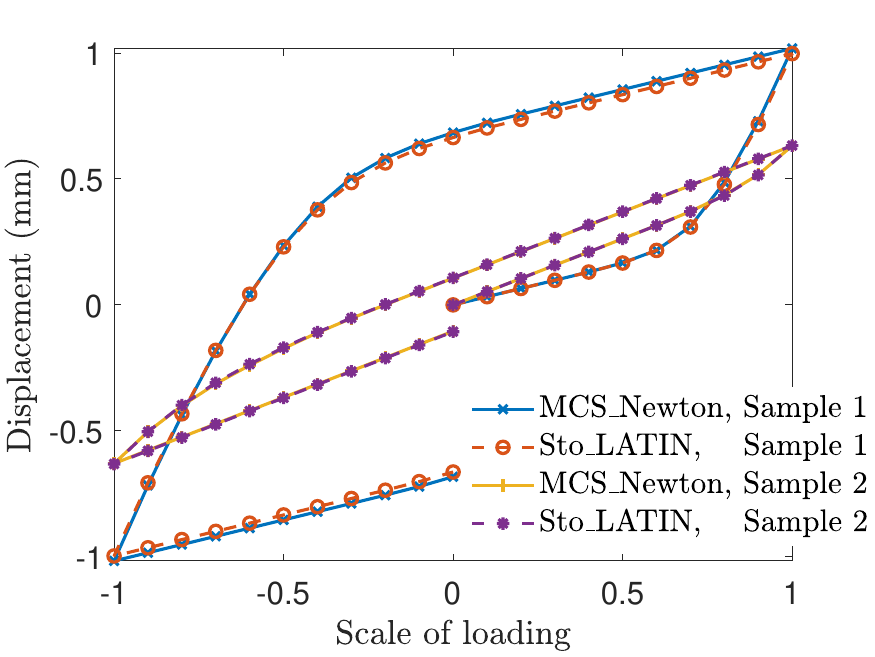}}
    \end{minipage}
    \hspace{0.12cm}
    \begin{minipage}{0.5\textwidth}
        \centering
        \subfloat[][Absolute errors of the displacements referring to the MCS reference solutions.]{\label{fig_e1_Dis_A_2_sam_Ab_Err} \includegraphics[width=1.0\textwidth]{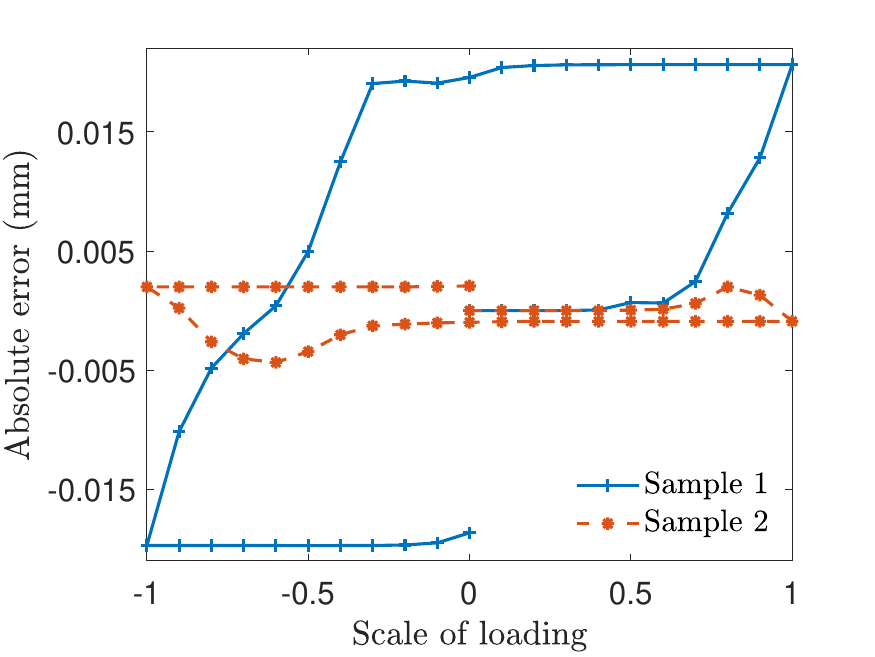}}
    \end{minipage}
    \caption{The stochastic displacements of the reference point $A$ in the $y$ direction of two sample realizations obtained by the Newton-based MCS and the stochastic LATIN iterations without and with update (left) and the absolute errors referring to the MCS reference displacements (right).} \label{fig_e1_Dis_A_2_sam_err}
\end{figure}

    The computational time of the stochastic LATIN method is 3.79 hours, including the iteration time of 0.43 hours and the update time of 3.36 hours for $1 \times 10^3$ random sample realizations.
    As a comparison, the MCS spends 60.79 hours, which is much higher than the proposed method.
    Further, two new sample realizations of random inputs, $\left\{ E\left( \theta^{\left( 1 \right)} \right), \sigma_{\rm Y}\left( \theta^{\left( 1 \right)} \right) \right\} = \left\{ 2.57 \times 10^5, 192.62 \right\}$ and $\left\{ E\left( \theta^{\left( 2 \right)} \right), \sigma_{\rm Y}\left( \theta^{\left( 2 \right)} \right) \right\} = \left\{ 1.59 \times 10^5, 291.40 \right\}$ that do not belong to the 20 random samples used to the iterative process, are used to test the generalization performance of the proposed method to new random sample inputs.
    The stochastic displacements and their absolute errors are depicted in \figref{fig_e1_Dis_A_2_sam_err} and corresponding $L_2$ errors are respectively $5.47 \times 10^{-4}$ and $2.40 \times 10^{-5}$, which demonstrates the good generalization of the stochastic LATIN method.
    In this sense, the proposed method also provides a high-accuracy stochastic reduce-order model of the stochastic problem under consideration.

\subsection{Example 2: Parameterized elastoplastic problem}
    In this example, we consider a parameterized elastoplastic problem and the Young's modulus and the yield stress are described using the following bounded interval parameters

    \vspace{-1.0em}
    \begin{equation*}
        E \in \left[ E_1, E_2 \right] = \left[ 1.90, 2.32 \right] \times 10^5~{\rm MPa}, \quad 
        \sigma_{\rm Y} \in \left[ \sigma_{{\rm Y},1}, \sigma_{{\rm Y},2} \right] = \left[ 220.50, 269.50 \right]~{\rm MPa}.
    \end{equation*}

\begin{figure}[ht]
    \centering
    \includegraphics[width=0.5\linewidth]{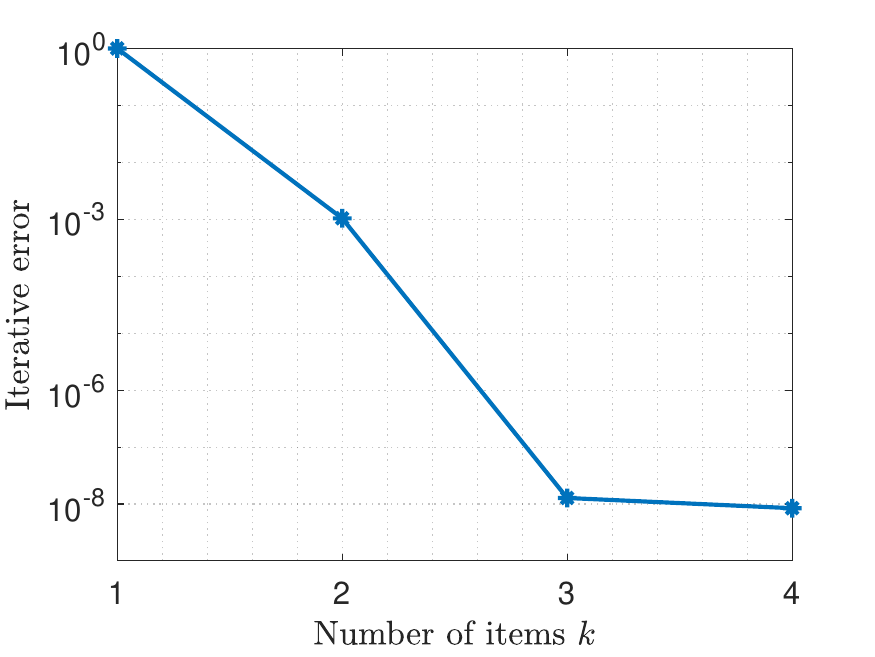}
    \caption{Iterative errors of diﬀerent retained terms.}
    \label{fig_e2_It_err}
\end{figure}

    As mentioned above, parameters on the given intervals are treated as uniform random variables in the proposed stochastic LATIN framework.
    The sample size $n_{s,1} = 20$ is again adopted for the iteration of Algorithm \ref{Alg_Sto_LATIN}.
    Iterative errors of different retained terms are shown in \figref{fig_e2_It_err} and only 4 terms are retained in this case.
    Since the ranges of parameter values are smaller than the ranges of the random variable values in Example \ref{subsection:Sto_Case}, fewer retained terms are sufficient to approximate the stochastic solution.
    Also, the linear global equation (\ref{eq:KdF}) is solved 11 times.
    The iteration time is 0.18 hours, which is much cheaper than Example \ref{subsection:Sto_Case} due to fewer triplets retained.

\begin{figure}[ht]
    \begin{minipage}{0.5\textwidth}
        \centering
        \subfloat[][The stochastic displacements obtained by the MCS and the stochastic LATIN iteration without and with update.]{\label{fig_e2_Dis_A_01} \includegraphics[width=1.0\textwidth]{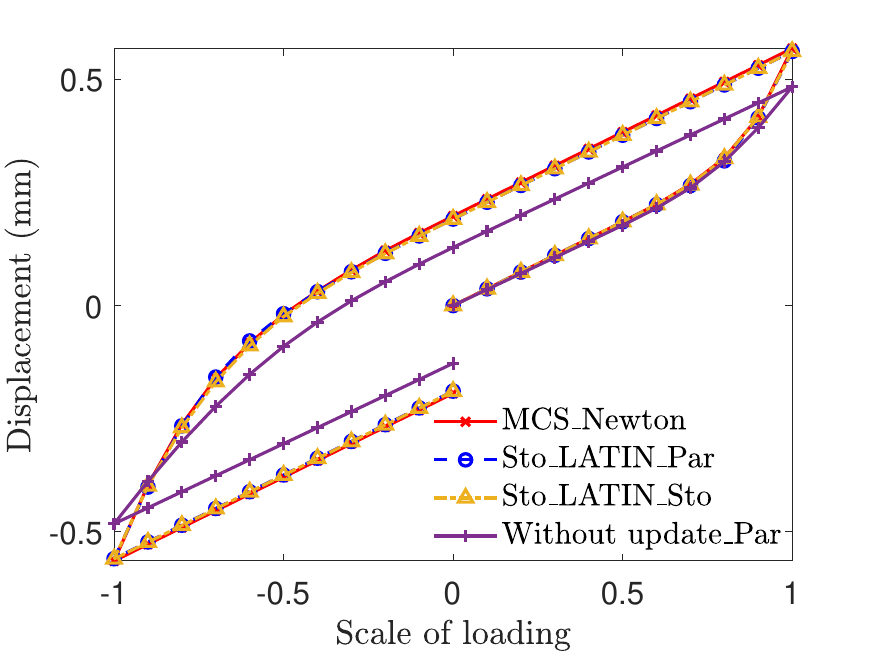}}
    \end{minipage}
    \hspace{0.12cm}
    \begin{minipage}{0.5\textwidth}
        \centering
        \subfloat[][Absolute errors of the displacements referring to the MCS reference solutions.]{\label{fig_e2_Dis_A_01_Ab_Err} \includegraphics[width=1.0\textwidth]{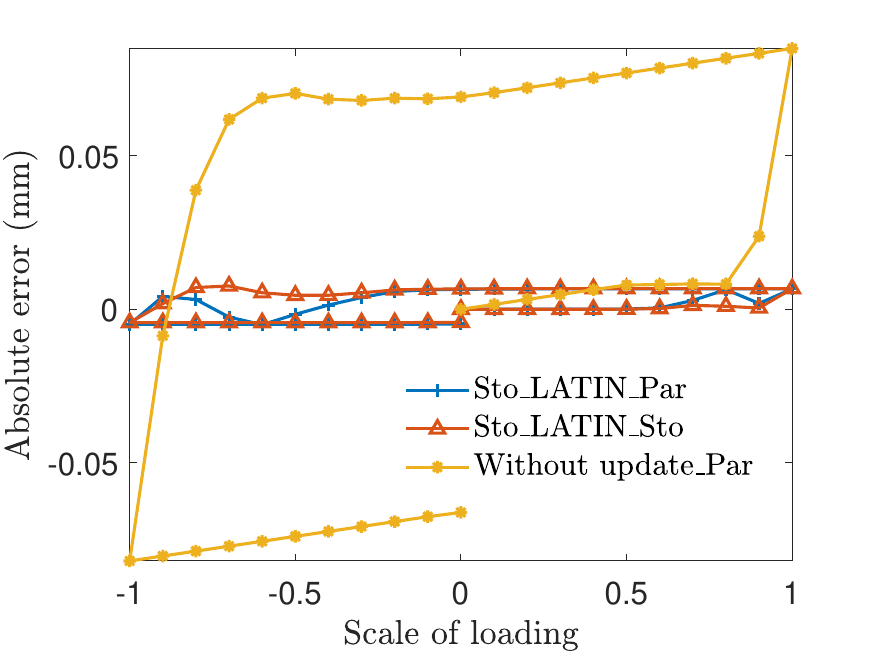}}
    \end{minipage}
    \caption{The stochastic displacements (left) of the point $A$ in the $y$ direction of the sample realization $\left\{ E^*, \sigma_{\rm Y}^* \right\} = \left\{ 2.26 \times 10^5, 258.73 \right\}$ obtained by the MCS and the stochastic LATIN method without and with update and the absolute errors referring to the Newton-based MCS reference displacement (right).} \label{fig_e2_Dis_A_01_err}
\end{figure}

\begin{figure}[!ht]
    \centering
    \includegraphics[width=1.0\linewidth]{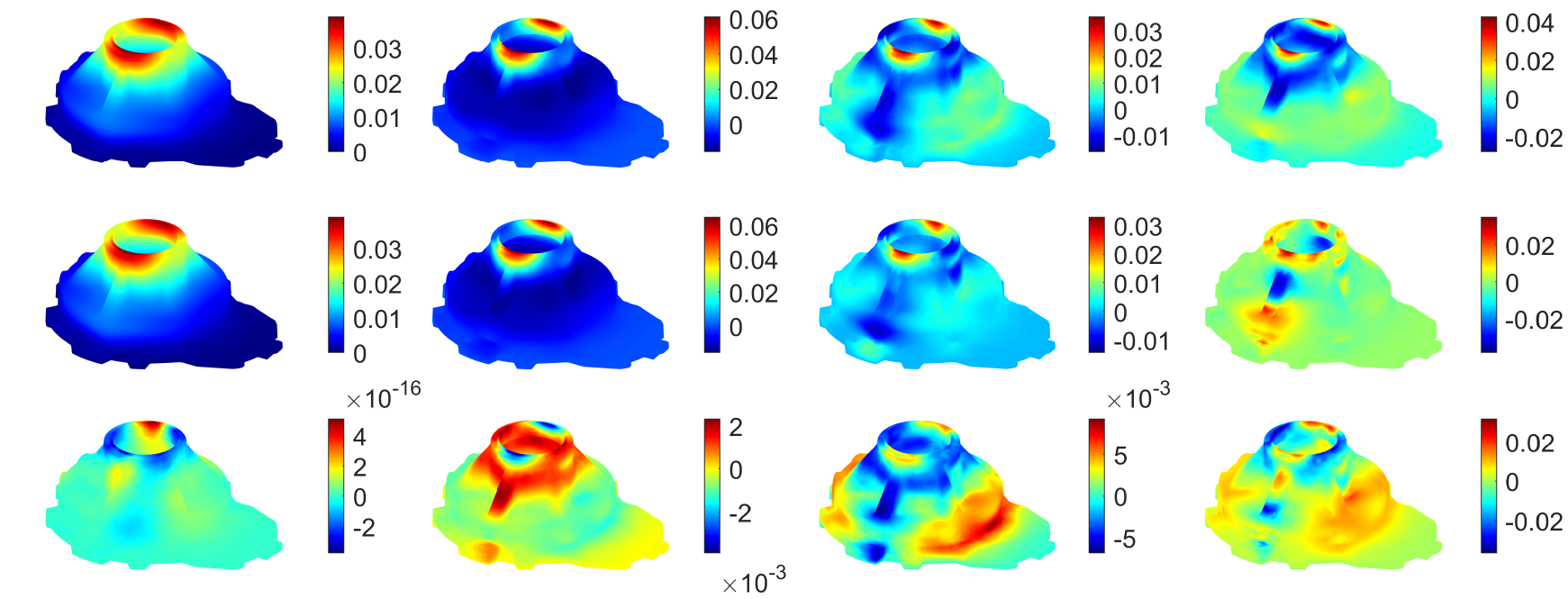}
    \caption{Comparison between the first four spatial functions $\left\{ {\bf d}_{i,y} \right\}_{i=1}^4$ in the $y$ direction calculated by the random inputs (i.e. the case in Example \ref{subsection:Sto_Case}) (the first line) and the parameter inputs (the second line) in this example, and their absolute errors (the third line).} 
    \label{fig_e2_D_Com}
\end{figure}

    We highlight that the ranges of parameter values are subsets of ranges of the random variable values in Example \ref{subsection:Sto_Case}.
    Thus, the spatial functions obtained in Example \ref{subsection:Sto_Case} can also be applied to the update stage of this example.
    To illustrate this and test the computational accuracy of the proposed method for parametric problems, the stochastic displacements of the point $A$ in the $y$ direction of the sample realization $\left\{ E^*, \sigma_{\rm Y}^* \right\} = \left\{ 2.26 \times 10^5, 258.73 \right\}$ are solved by using the Newton iteration-based MCS, the stochastic LATIN method without and with update under parametric inputs (given in this example), and the stochastic LATIN method under random inputs (given in Example \ref{subsection:Sto_Case}), respectively.
    Their comparisons are found in \figref{fig_e2_Dis_A_01_err}, and the $L_2$ errors of the stochastic LATIN method without and with update under parametric inputs, and the stochastic LATIN method under random inputs are $3.75 \times 10^{-2}$, $2.22 \times 10^{-4}$ and $2.34 \times 10^{-4}$, respectively, which indicates that the recalculation process is necessary to improve the accuracy of the stochastic solution for both random and parametric inputs, and also verifies the effectiveness of the reuse of spatial functions from Example \ref{subsection:Sto_Case}.
    Further, the computational times of the stochastic displacement solved by the MCS, and the update stages of stochastic LATIN method under parametric and random inputs are 189.84s, 8.15s and 8.69s, respectively.
    The computational cost of the proposed method is greatly reduced compared with the MCS.
    The stochastic LATIN method under parametric inputs is slightly cheaper than that under random inputs due to fewer triplets involved.
    The first four spatial functions in the $y$ direction obtained in Example \ref{subsection:Sto_Case} and this example are depicted in the first and second lines of \figref{fig_e2_D_Com}, and their absolute errors are found in the third line of \figref{fig_e2_D_Com}.
    It is seen that the first three dominant spatial functions have very similar modes.
    In this sense, the proposed method is less sensitive to the distributions of inputs.
    The random sample realizations used in the iterative process can hence be generated without careful design.

\subsection{Example 3: Mixed stochastic and parameterized elastoplastic problem}
    In this example, we further consider an elastoplastic problem with mixed random and parametric inputs.
    Specifically, the Young's modulus and the yield stress are respectively modeled as a Gaussian random variable and an interval parameter

    \vspace{-1.0em}
    \begin{equation*}
        E\left( \theta \right) \sim {\cal N}\left( m_{E}, \chi_E^2 \right), \quad 
        \sigma_{\rm Y} \in \left[ \sigma_{{\rm Y},1}, \sigma_{{\rm Y},2} \right] = \left[ 220.50, 269.50 \right]~{\rm MPa},
    \end{equation*}
    where the mean value $m_E = 2.11 \times 10^5$~MPa and the standard deviation $\chi_E = 0.1 m_E$.

\begin{figure}[ht]
    \begin{minipage}{0.5\textwidth}
        \centering
        \subfloat[][PDFs of the time step $t_{11}$.]{\label{fig_e3_PDF_t_11} \includegraphics[width=1.0\textwidth]{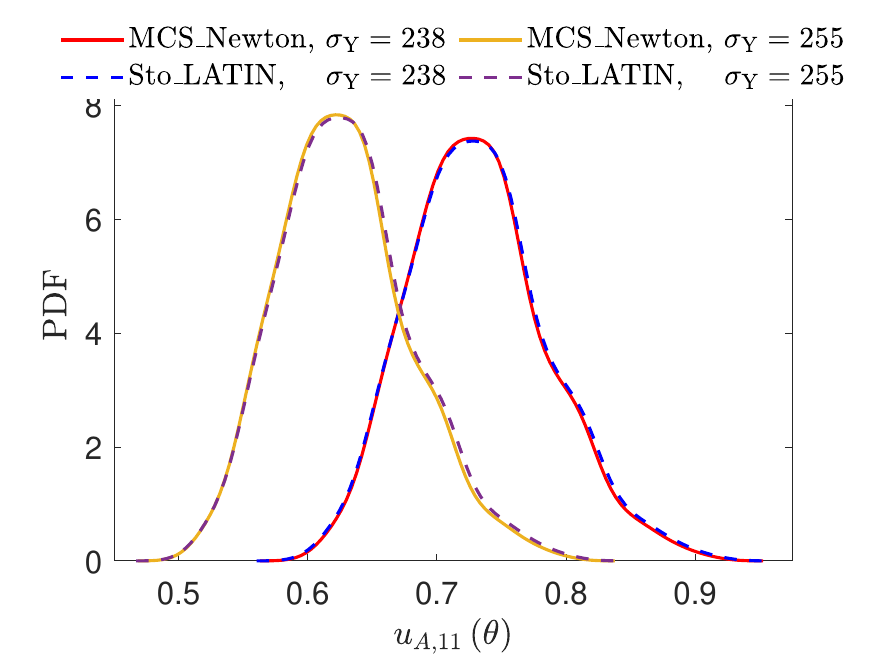}}
    \end{minipage}
    \hspace{0.12cm}
    \begin{minipage}{0.5\textwidth}
        \centering
        \subfloat[][PDFs of the time step $t_{41}$.]{\label{fig_e3_PDF_t_41} \includegraphics[width=1.0\textwidth]{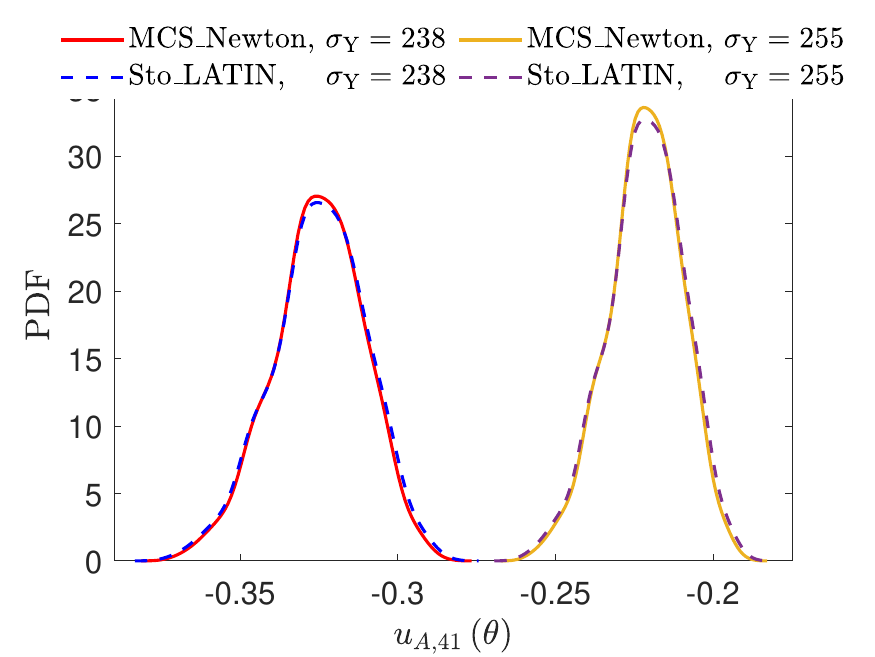}}
    \end{minipage}
    \caption{PDFs of the stochastic displacement of the point $A$ in the $y$ direction at time steps $t=11$ (left) and $t=41$ (right) obtained by the MCS and the proposed stochastic LATIN method.} \label{fig_e3_PDF_t_11_41}
\end{figure}

    We do not solve this example by the iteration in Algorithm \ref{Alg_Sto_LATIN}, but only perform the update stage by reusing the spatial functions obtained in Example \ref{subsection:Sto_Case}.
    The stochastic solutions are solved under the random Young's modulus mixed with two sample realizations of the yield stress, given by $\left\{ E\left( \theta \right), \sigma_{\rm Y} = 238 \right\}$ and $\left\{ E\left( \theta \right), \sigma_{\rm Y} = 255 \right\}$. 
    PDFs of the stochastic displacement of the point $A$ in the $y$ direction at time steps 11 and 41 obtained by the MCS and the stochastic LATIN method are shown in \figref{fig_e3_PDF_t_11} and \figref{fig_e3_PDF_t_41}, respectively, which demonstrates that the proposed method still has very high accuracy for mixed random and parametric inputs.
    For both time steps $t_{11}$ and $t_{41}$, there are larger stochastic displacements in the case of $\sigma_{\rm Y} = 238$ since larger plastic deformation occurs under the smaller yield stress $\sigma_{\rm Y} = 238$ compared to the case $\sigma_{\rm Y} = 255$.
    Computational times of the two cases are listed in Table \ref{tab_e3_Time}, where the Sto\_LATIN time is the sum of the iteration time and the update time.
    The proposed method is much more efficient than the MCS for both cases.
    The two cases share the spatial functions from Example \ref{subsection:Sto_Case}, and thus have the same iteration times.
    The case $\sigma_{\rm Y} = 238$ takes a bit more computational time because more effort is put into solving the larger nonlinear plastic displacements.

\begin{table}[ht]
    \centering
    \caption{Computational times (hour) of two sample realizations $\sigma_{\rm Y} = 238$~MPa and 255~MPa.}\label{tab_e3_Time}
    \begin{tabular}{ccccc}
        \toprule
        $\sigma_{\rm Y}$ & Iteration time & Update time & Sto\_LATIN time & MCS\_Newton time \\
        \midrule
        238  & 0.43 & 3.43 & 3.86 & 61.43 \\
        255  & 0.43 & 3.29 & 3.72 & 58.50 \\
        \bottomrule
    \end{tabular}
\end{table}

\subsection{Example 4: High-dimensional stochastic elastoplastic problem with spatial variability}

    In this example, we consider a stochastic elastoplastic problem with high-dimensional random field inputs. 
    The Young's modulus $E\left(x,y,z,\theta\right)$ and the yield stress $\sigma_{\rm Y}\left(x,y,z,\theta\right)$ are modeled as two random fields with means $E_0 = 2.11 \times 10^5$~MPa and $\sigma_{{\rm Y},0} = 245$~MPa, respectively.
    Their covariance functions are given by a unified form

    \vspace{-1.0em}
    \begin{equation} \label{eq:ex4_Cov}
        {\rm Cov}_{\Box}\left( x_1,y_1,z_1; x_2,y_2,z_2 \right) = \frac{1}{3} \chi_{\Box}^2 \exp\left( -\frac{\left| x_1-x_2 \right|}{l_x} -\frac{\left| y_1-y_2 \right|}{l_y} -\frac{\left| z_1-z_2 \right|}{l_z} \right),
    \end{equation}
    where $\chi_{\Box}$ is 0.1 times the mean value, that is, $\chi_E = 0.1 E_0$ and $\chi_{\sigma_{\rm Y}} = 0.1 \sigma_{{\rm Y},0}$.
    The correlation lengths $l_x$, $l_y$ and $l_z$ are given by the maximum length of the physical model along the $x$, $y$ and $z$ directions, respectively.
    We let the above random fields approximate by the following series expansion 

    \vspace{-1.0em}
    \begin{equation} \label{eq:ex4_E}
        E\left(x,y,z,\theta\right) = E_0 + \sum\limits_{i=1}^{r_E} \xi_i\left( \theta \right) \sqrt{\kappa_{E,i}} E_i\left( x,y,z \right), \quad
        \sigma_{\rm Y}\left(x,y,z,\theta\right) = \sigma_{{\rm Y},0} + \sum\limits_{i=1}^{r_{\sigma_{\rm Y}}} \eta_i\left( \theta \right) \sqrt{\kappa_{\sigma_{\rm Y},i}} \sigma_{{\rm Y},i}\left( x,y,z \right),
    \end{equation}
    where $\left\{ \xi_i\left( \theta \right) \right\}_{i=1}^{r_E}$ and $\left\{ \eta_i\left( \theta \right) \right\}_{i=1}^{r_{\sigma_{\rm Y}}}$ are mutually independent uniform random variables on $\left[ -\sqrt{3}, \sqrt{3} \right]$.
    $\left\{ \kappa_{{\Box},i} {\Box}_i\left( x,y,z \right) \right\}_{i=1}^{r_{\Box}}$ are eigenvalues and eigenvectors of the covariance function and solved by the following Fredholm integral equation of the second kind

\begin{figure}[ht]
    \centering
    \includegraphics[width=0.5\linewidth]{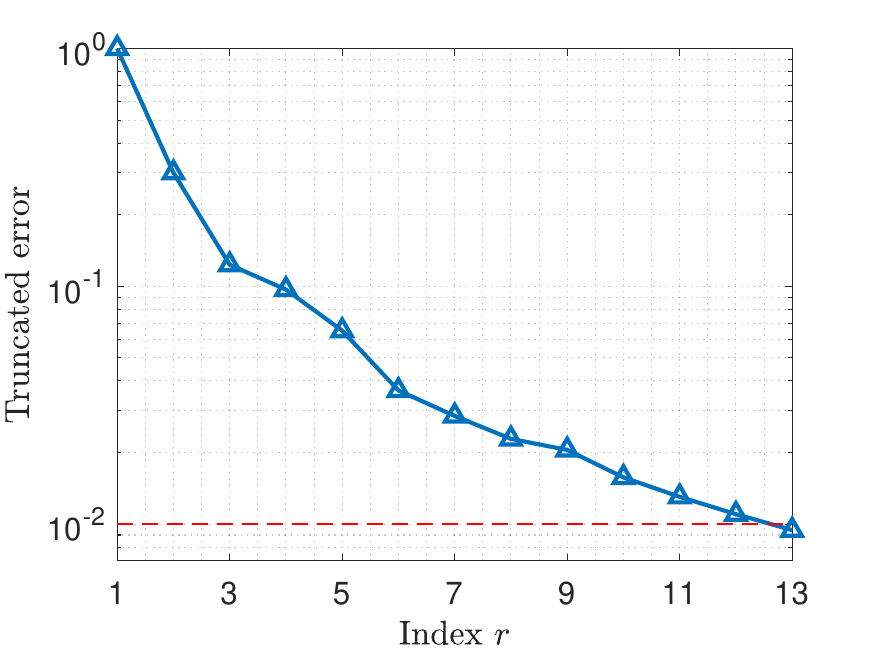}
    \caption{Truncated errors of diﬀerent truncation numbers $r$.}
    \label{fig_e4_Trun_err}
\end{figure}

\begin{figure}[ht]
    \centering
    \includegraphics[width=0.9\linewidth]{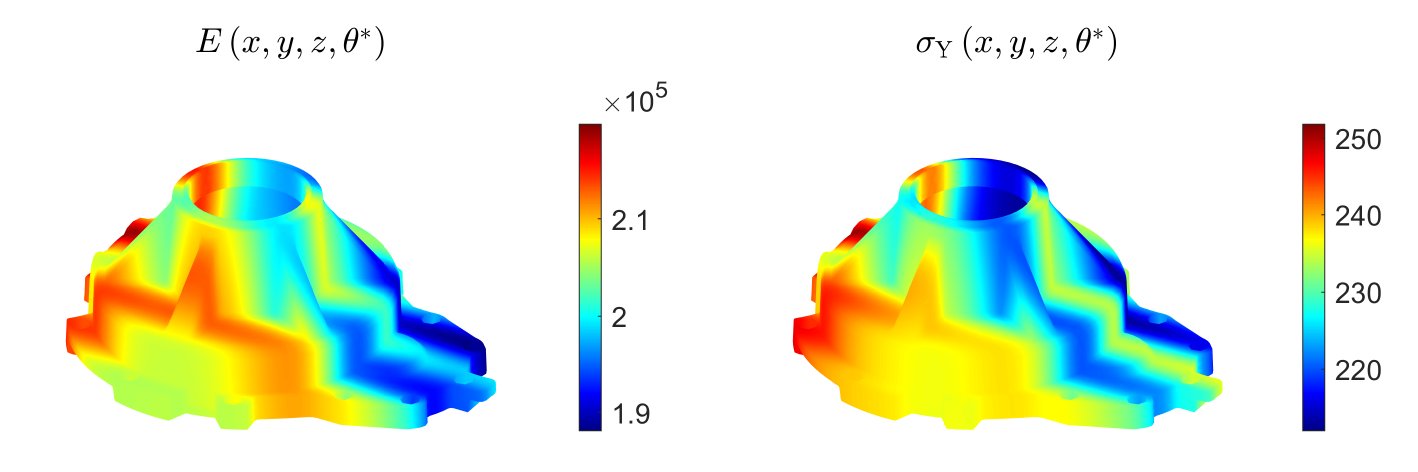}
    \caption{Sample realizations of the random fields: the Young's modulus $E\left( x,y,z, \theta^* \right)$ (left) and the yield stress $\sigma_{\rm Y}\left( x,y,z, \theta^* \right)$ (right).}
    \label{fig_e4_Sam_EY}
\end{figure}

    \vspace{-1.0em}
    \begin{equation}
        \int_{\Omega} {\rm Cov}_{\Box}\left( x_1,y_1,z_1; x_2,y_2,z_2 \right) {\Box}_i\left( x_1,y_1,z_1 \right) {\rm d}x_1{\rm d}y_1{\rm d}z_1 = \kappa_{{\Box},i} {\Box}_i\left( x_2,y_2,z_2 \right),
    \end{equation}
    which can be solved by existing numerical solvers \cite{saad2011numerical}.
    The number $r$ of truncation terms is determined by the truncated error ${{\kappa _{\Box,r}}} \mathord{\left/
    {\vphantom {{{\kappa _{\Box,r}}} {\sum\nolimits_{i = 1}^r {{\kappa _{\Box,i}}} }}} \right.
    \kern-\nulldelimiterspace} {\sum\nolimits_{i = 1}^r {{\kappa _{\Box,i}}} } \le 1 \times 10^{- 2}$.
    Truncated errors of diﬀerent truncation numbers $r$ are shown in \figref{fig_e4_Trun_err}.
    For this example, $r=13$ is adopted to achieve the given precision.
    Therefore, a total of 26 stochastic dimensions are involved in the stochastic elastoplastic analysis.
    Sample realizations of the Young's modulus and the yield stress are depicted in \figref{fig_e4_Sam_EY}, which exhibits high spatial variability.

\begin{figure}[ht]
    \centering
    \includegraphics[width=0.5\linewidth]{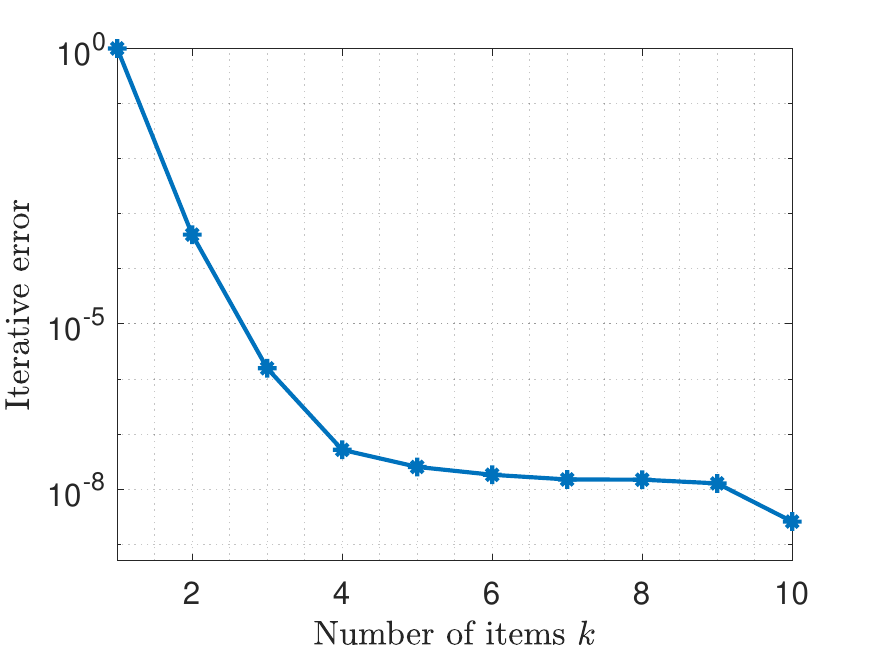}
    \caption{Iterative errors of diﬀerent retained terms.}
    \label{fig_e4_It_err}
\end{figure}

\begin{figure}[!ht]
    \begin{minipage}{0.5\textwidth}
        \centering
        \subfloat[][The stochastic displacements obtained by the MCS and the stochastic LATIN iteration without and with update.]{\label{fig_e4_Dis_A_01} \includegraphics[width=1.0\textwidth]{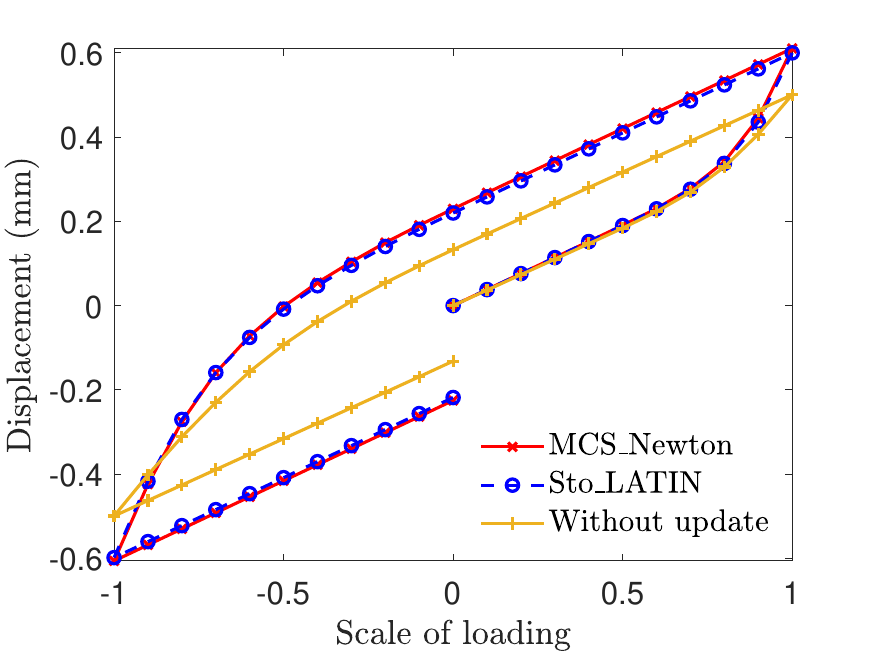}}
    \end{minipage}
    \hspace{0.12cm}
    \begin{minipage}{0.5\textwidth}
        \centering
        \subfloat[][Absolute errors of the displacements referring to the MCS reference solutions.]{\label{fig_e4_Dis_A_01_Ab_Err} \includegraphics[width=1.0\textwidth]{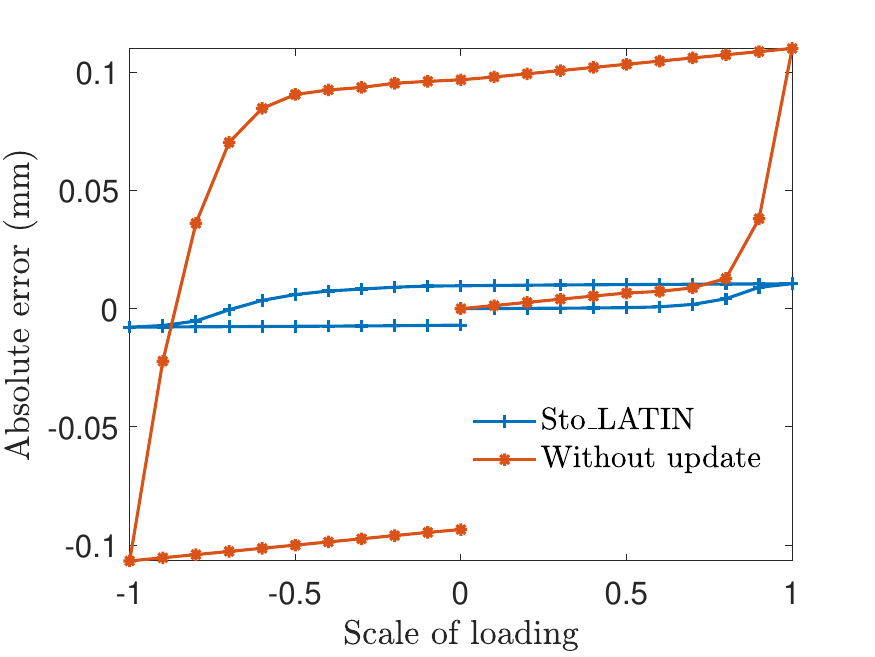}}
    \end{minipage}
    \caption{The stochastic displacements of the point $A$ in the $y$ direction of the sample realization $\left\{ E^*, \sigma_{\rm Y}^* \right\} = \left\{ 2.26 \times 10^5, 258.73 \right\}$ obtained by the MCS and the stochastic LATIN iterations without and with update (left) and the absolute errors referring to the Newton-based MCS reference displacement (right).} \label{fig_e4_Dis_A_01_err}
\end{figure}

    In this case, the sample size $n_{s,1} = 260$ is set for the iterative process of Algorithm \ref{Alg_Sto_LATIN}.
    Iterative errors of different retained terms are plotted in \figref{fig_e4_It_err}, and 10 triplets are retained to meet the specified convergence error.
    Compared to the previous examples, more terms are required to approximate the stochastic solution due to the spatial variability of random inputs.
    Also, 52 linear global deterministic equations are solved with an iteration time of 2.81 hours, which is more expensive than the above problems with random inputs without spatial variability.
    Further, the computational accuracy for this case is checked using the stochastic displacement corresponding to the sample realization shown in \figref{fig_e4_Sam_EY}.
    Specifically, the stochastic displacements of the point $A$ in the $y$ direction are solved by using the MCS, and the stochastic LATIN method without and with update, respectively.
    The stochastic displacement comparisons can be found in \figref{fig_e4_Dis_A_01_err}, and the corresponding $L_2$ errors of the stochastic LATIN method without and with update are $5.70 \times 10^{-2}$ and $4.51 \times 10^{-4}$, respectively.
    The computational times of the stochastic displacement solved by the MCS and the update stage of stochastic LATIN method are 262.73s and 8.31s, respectively.
    Therefore, both good accuracy and high efficiency are achieved for the high-dimensional stochastic problem with spatial variability.

\begin{figure}[!ht]
    \centering
    \includegraphics[width=0.5\linewidth]{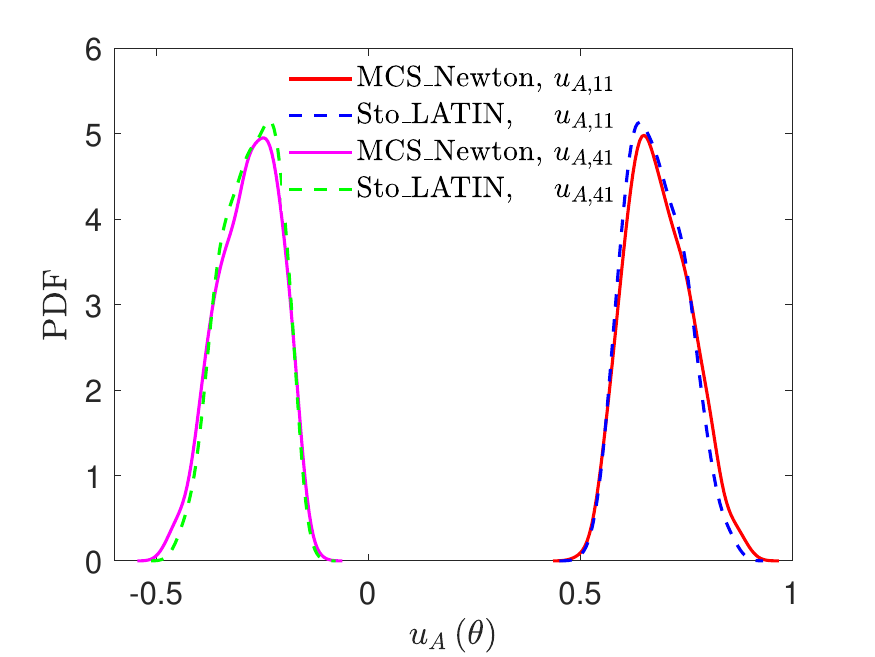}
    \caption{PDFs of stochastic displacement of the point $A$ in the y direction at time steps $t=11$ and $t=41$ obtained by the MCS and the stochastic LATIN method.}
    \label{fig_e4_PDF_A_11_41}
\end{figure}

\begin{figure}[ht]
    \begin{minipage}{0.5\textwidth}
        \centering
        \subfloat[][Iterative errors of diﬀerent retained terms.]{\label{fig_e4_Restart_err} \includegraphics[width=1.0\textwidth]{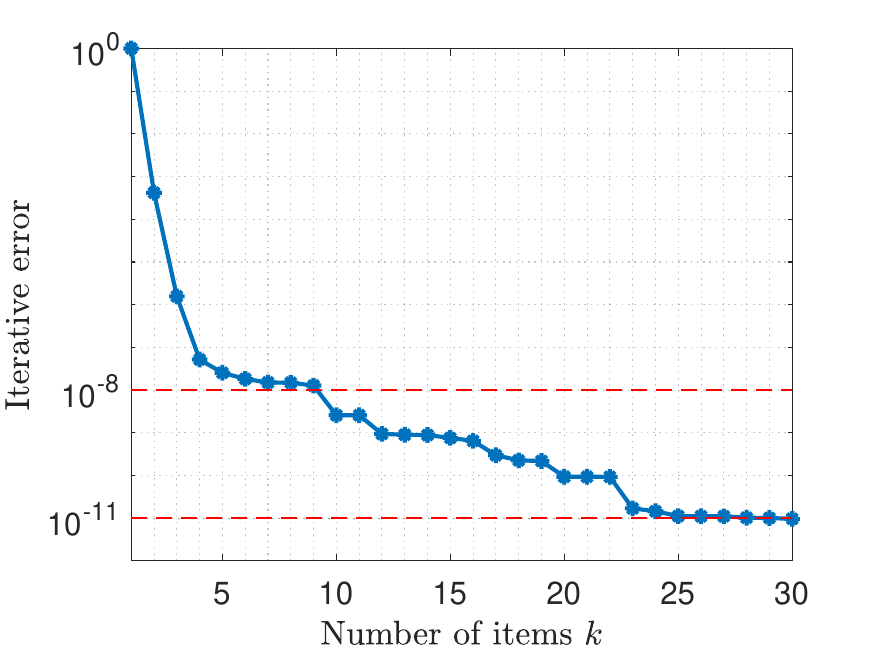}}
    \end{minipage}
    \hspace{0.12cm}
    \begin{minipage}{0.5\textwidth}
        \centering
        \subfloat[][PDFs of the time steps $t_{11}$ and $t_{41}$.]{\label{fig_e4_Restart_PDF} \includegraphics[width=1.0\textwidth]{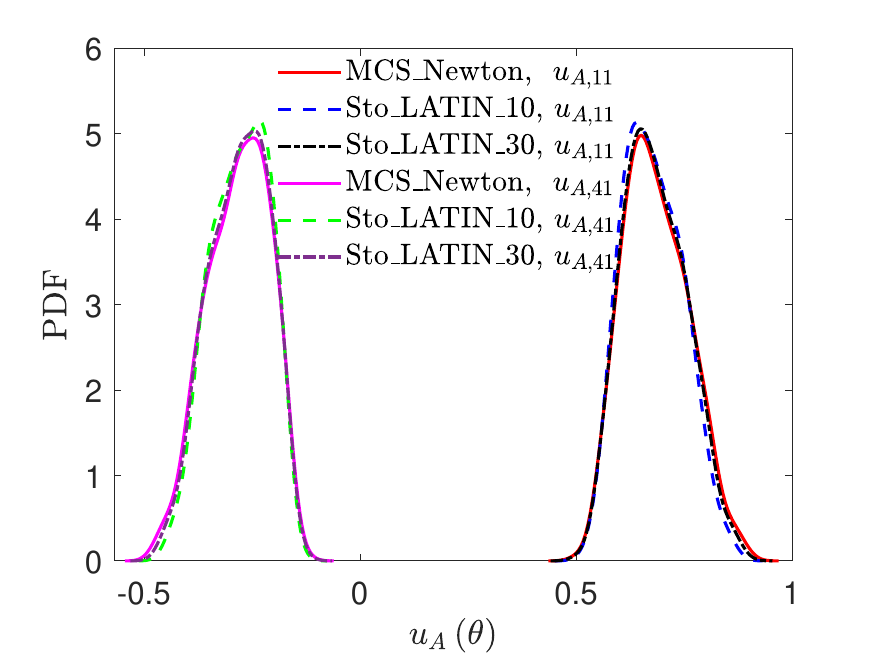}}
    \end{minipage}
    \caption{Restart Algorithm \ref{Alg_Sto_LATIN} to reach the convergence error $1 \times 10^{-11}$: Iterative errors of different retained terms (left), and PDFs of the displacement the point $A$ in the y direction at time steps $t=11$ and $t=41$ obtained by the MCS and the stochastic LATIN method (right).} \label{fig_e4_Restart}
\end{figure}

    PDFs of stochastic displacement of the point $A$ in the $y$ direction at the time steps 11 and 41 obtained by $1 \times 10^3$ MCS and the proposed stochastic LATIN method are compared in \figref{fig_e4_PDF_A_11_41}.
    The proposed method is still consistent with the MCS, but the approximation accuracy is not as good as that in the previous examples. 
    Due to the spatial variability of the stochastic solution, more triplets are required to achieve a higher-accuracy approximation.
    The proposed method only takes 6.49 hours (including the iteration time of 2.81 hours and the update time of 3.68 hours), which is much cheaper than the cost of MCS of 62.08 hours.
    Further, we highlight that the proposed stochastic LATIN method is easily restarted to achieve higher-accuracy stochastic solutions.
    More triplets can be continuously solved by Algorithm \ref{Alg_Sto_LATIN}, where the obtained triplets $\left\{ \lambda_i\left( \theta \right), g_i\left( t \right), {\bm d}_i\left( {\bf x} \right) \right\}_{i=1}^k$ are treated as the initialization of Algorithm \ref{Alg_Sto_LATIN}.
    To illustrate this, we restart Algorithm \ref{Alg_Sto_LATIN} to achieve a more accurate approximation of the above PDFs.
    To this end, we reset the convergence error of stochastic solution is $\epsilon _{\bf u} = 1 \times 10^{-11}$.
    Iterative errors of different retained terms and the improved PDFs obtained by the MCS and the stochastic LATIN method are shown in \figref{fig_e4_Restart_err} and \figref{fig_e4_Restart_PDF}, respectively.
    20 additional terms are solved based on the obtained 10 triplets, and a total of 30 triplets are retained to approximate the stochastic solution.
    PDFs are thus greatly improved and agree very well with the MCS reference PDFs.

\section{Conclusions and discussions} \label{section:Conclusion}
    This paper presented an effective stochastic LATIN framework for solving elastoplastic problems with random and/or parametric inputs.
    A decoupling representation is proposed to approximate stochastic solutions using a set of triplets of spatial functions, temporal functions and random variables.
    A stochastic LATIN iteration is developed to calculate each triplet alternatively.
    The high efficiency of the proposed method is reached since only linear global equations are solved, and good accuracy is achieved by controlling the number of retaining terms, which has been verified by 3D numerical examples with random variables, parameters and high-dimensional random field inputs.
    Although only the stochastic elastoplastic problem is considered in this paper, the proposed framework can be extended to more general nonlinear stochastic problems.
    It also provides a novel and potential method for uncertainty quantification in science and engineering.
    However, a preselection of the sample size for the iterative process of the proposed stochastic LATIN framework remains an open problem, which will be investigated in depth in subsequent research.

\section*{Acknowledgments}
    The authors are grateful to the Alexander von Humboldt Foundation and the International Research Training Group 2657 (IRTG 2657) funded by the German Research Foundation (DFG) (Grant number 433082294).

\nocite{*}
\bibliography{References}

\end{document}